\documentclass[11pt]{amsart}
\usepackage{ amsthm, amssymb, amsfonts, graphicx}
\usepackage{comment}
\usepackage{inputenc}
\usepackage{mathtools}
\usepackage[margin= 1.0 in]{geometry}
\usepackage{xcolor}
\usepackage{enumitem}

\newtheorem{thm}{Theorem}[section]
\newtheorem{lem}[thm]{Lemma}
\newtheorem{prop}[thm]{Proposition}
\newtheorem{cor}[thm]{Corollary}

\usepackage{tikz-cd}
\usepackage[
backend=biber,
style=alphabetic,
]{biblatex}

\bibliography{ADref.bib}

\theoremstyle{definition}
\newtheorem{defn}{Definition}[section]

\theoremstyle{remark}
\newtheorem{rem}{Remark}

\numberwithin{equation}{section}

\def\C{{\mathbb C}}

\def\Q{{\mathbb Q}}
\def\P{{\mathbb P}}
\def\Q{{\mathbb Q}}
\def\R{{\mathbb R}}

\def\0{{\mathbf 0}}
\def\1{{\mathbf 1}}

\def\w{{\mathbf w}}

\def\Jcal{{\mathcal J}}
\def\Kcal{{\mathcal K}}

\def\Gal{\mathrm{Gal}}
\def\Gauss{\mathrm{Gauss}}

\def\PGL{\mathrm{PGL}}

\def\Res{\mathrm{Res}}
\def\Rat{\mathrm{Rat}}

\def\supp{\mathrm{supp}}

\def\BDV{\mathrm{BDV}}
\def\Res{\mathrm{Res}}

\def\Re{\mathrm{Re}\,}
\def\<{\langle}
\def\>{\rangle}

\def\nubar{\overline{\nu}}

\def\Pberk{ \mathbf{P}_{v}^1 }
\def\Aberk{ \mathbf{A}_{ v}^1 }
\def\Ar{\mathrm{Ar}}
\def\AZ{\mathrm{AZ}}
\def\std{\mathrm{std}}
\def\Aad{\mathcal{A}_K}

\def\Diag{ \mathrm{Diag}}

\DeclareMathOperator{\qif}{\;\; \text{if} \;\;}

\DeclareMathOperator{\qas}{\;\; \text{as} \;\;}
\DeclareMathOperator{\qand}{\;\;\text{and} \;\;}

\title[An Inner Product on Adelic Measures]{%
An Inner Product on Adelic Measures \\ With Applications to the Arakelov--Zhang Pairing
}
\author{Peter J. Oberly}
\begin{document}

\maketitle


\begin{abstract}
We define an inner product on a vector space of adelic measures over a number field. We find that the norm induced by this inner product governs weak convergence at each place of $K$. The canonical adelic measure associated to a rational map is in this vector space, and the square of the norm of the difference of two such adelic measures is the Arakelov--Zhang pairing from arithmetic dynamics. We find that the norm of a canonical adelic measure associated to a rational map is commensurate with a height on the space of rational functions with fixed degree. As a consequence, we show that the Arakelov--Zhang pairing of two rational maps $f$ and $g$ is, when holding $g$ fixed, commensurate with the height of $f$. 
\end{abstract}

\section{Introduction}\label{sec: intro}

Let $K$ be a number field and let $M_K$ be the set of places of $K$. By an \textit{adelic measure}, we mean a sequence $\mu = (\mu_v)_{v \in M_K}$ of suitably regular complex-valued measures indexed by the places of $K$. At the finite places of $K$, each measure $\mu_v$ is supported on the Berkovich projective line $\Pberk$ and $\mu_v$ is supported on the Riemann sphere at the infinite places. We also impose two global conditions on adelic measures, namely that the collection of measures $(\mu_v)_{v \in M_K}$ has constant total mass, and that $\mu_v$ is trivial for all but finitely many places; precise definitions are deferred to Section \ref{sec: adelic measures}. The collection of such adelic measures forms complex vector space, which we denote by $\Aad$.\\

The principal example of an adelic measure is the canonical adelic measure $\mu_f = (\mu_{f, v})_{v \in M_K}$ of a rational map $f$ defined over $K$. This adelic measure encodes a good amount of information on the $v$-adic Julia sets of $f$, and so defines an important dynamical invariant. For example, the preperiodic points of $f$ distribute according to $\mu_{f, v}$ at each place. When $v$ is an Archimedean place, $\mu_{f, v}$ is the measure investigated for polynomials by Brolin \cite{brolin1965invariant},  and for rational maps by Ljubich \cite{ljubich1983entropy} and by Freire, Lopes, and Ma\~n{\'e}, \cite{freire1983invariant}. For non-Archimedean places, this measure was introduced independently by several authors, including Baker--Rumely \cite{baker2006equidistribution}, Chambert-Loir \cite{chambert2006mesures}, and Favre--Rivera-Letelier \cite{favre2006equidistribution}. \\

 We equip $\Aad$ with the following pairing. Fix an adelic probability measure $\lambda = (\lambda_v)_{v \in M_K}$. At each place $v$ we let $g_{\lambda, v}(x,y)$ be an Arakelov--Green's function of $\lambda_v$. (See Section 2 for the definition). We choose a normalization of $g_{\lambda_v}$ so that
\[\int_{\Pberk}\int_{\Pberk} g_{\lambda, v}(x,y) \,d\lambda_v(x)\,\lambda_v(y) =: t_v \geq 0  \]
and require $t_v = 0$ for all but finitely many $v$, and that  $\sum_{v \in M_K} r_v t_v =: t > 0$, where $r_v = [K_v:\Q_v]/[K:\Q]$ is the ratio of local to global degrees. Inspired by the work of Favre--Rivera-Letelier and of Baker-Rumely (discussed below), given two adelic measures $\mu, \nu \in \Aad$ we define 

\[  \< \mu, \nu\>_{\lambda, t} = \sum_{v \in M_K} r_v \int_{\Pberk} \int_{\Pberk} g_{\lambda, v}(x,y) \, d\mu(x)\, d\overline{\nu_v}(y). \]

 For example, take  the \textit{Arakelov adelic measure}, $\lambda_\Ar = (\lambda_{\Ar,v})_{v \in M_K}$. At Archimedean places $v$, $\lambda_{\Ar, v}$ is the spherical measure, and is trivial at finite places (see (\ref{eq: definition of arakelov measure})). An Arakelov--Green's function for $\lambda_{\Ar,v}$ which satisfies the above conditions is $g_{\lambda_{\Ar},v}(x,y) = \log ||x, y||_v^{-1}$, where $||\cdot,\cdot||_v$ is a continuous extension of the projective metric on $\P^1(\C_v)$ to $\Pberk$. (This is called the \textit{spherical kernel} in the book by Baker and Rumely \cite{baker2010potential}). It is not hard to check that $\int_{\Pberk}\int_{\Pberk} \log ||x, y||_v^{-1}\, d\lambda_{\Ar, v}(x)\, d\lambda_{\Ar}(y) = 1/2$ for Archimedean $v$ and is $0$ for non-Archimedean $v$; thus $\sum_{v \in M_K} r_v t_v = 1/2$ and so 
\[ \< \mu, \nu\>_{\lambda_{\Ar}, 1/2} = \sum_{v \in M_K} r_v \int_{\Pberk}\int_{\Pberk} \log \frac{1}{||x, y||_v} \, d\mu_v(x)\, d\overline{\nu_v}(y). \]
 
Our first theorem is that this pairing is an inner product, and that, in most situations of interest, convergence in the norm $\|\cdot\|_{\lambda,t}$ induced by this inner product gives weak convergence at each place. 

\begin{thm}\label{thm: inner product into}
Fix an adelic probability measure $\lambda = (\lambda_v)_{v \in M_K}$ and a normalization of Arakelov--Green's functions as above. Then map $(\mu, \nu) \mapsto\<\mu,\nu\>_{\lambda,t}$ is an inner product on $\Aad$. Furthermore, if $\mu_n, \mu \in \Aad$ with $\sup_n |\mu_{n,v}|(\Pberk) < \infty$ for all $v \in M_K$, then $\|\mu_n - \mu\|_{\lambda,t} \to 0$ implies $\mu_{n, v} \to \mu_v$ weakly at each place.
\end{thm}

Similar potential theoretic pairings are used by both Favre--Rivera-Letelier and Baker--Rumely in connection with their respective proofs of the adelic equidistribution theorem \cite{favre2006equidistribution},\cite{baker2006equidistribution}. The novelty of Theorem \ref{thm: inner product into} is that $\< \cdot, \cdot\>_{\lambda, t}$ is positive definite for adelic measures of arbitrary mass. In particular, the pairing $\< \cdot, \cdot\>_{\lambda, t}$ should perhaps be viewed as a modest extension of the potential theoretic pairings investigated by Favre--Rivera-Letelier and Baker--Rumely.  The following comparisons between the height of a rational map $f$ and the norm $\|\mu_f\|^2_{\lambda,t}$ of the canonical adelic measure of $f$ are the primary new results of this article.

\begin{thm}\label{thm: heigh on monic poly}
The map $f \mapsto \frac{d}{2}\|\mu_f\|_{\lambda,t}^2$ is a Weil height on the space of monic polynomials of degree $d$ defined over $K$. 
\end{thm}

The proof of Theorem \ref{thm: heigh on monic poly} relies on studying the $v$-adic Julia sets of $f$ at each place, and uses escape rate arguments to bound $\frac{d}{2}\|\mu_f\|^2_{\lambda, t}$ in terms of the height of the coefficients of $f$. Such  arguments do not apply to rational maps. However, our main theorem shows that $ f \mapsto \frac{d}{2}\|\mu_f\|_{\lambda, t}^2$  is at least commensurate with a height on $\Rat_d(K)$, the space of rational maps with degree $d$ defined over $K$. 

\begin{thm}\label{thm: height on polynomials}
Let $h$ be a Weil height on $\Rat_d(K)$. There are positive constants $c_1, c_2, c_3, c_4$  depending only on $h,\lambda,t$, and $d$ so that $c_1 h(f) - c_2\leq \frac{d}{2}\|\mu_f\|_{\lambda,t}^2 \leq c_3 h(f) + c_4$ for all $f \in \Rat_d(K)$. Moreover $c_1, c_3 \to 1$ as $d \to \infty$ and $c_2, c_4$ grow linearly in $d$.  
\end{thm}
  
These constants are given explicitly in the case of the Arakelov height of the coefficients of $f$ and $\|\mu_f\|^2_{\lambda_\Ar,1/2}$ (see (\ref{eq: explicit bounds arakelov height})). The main idea of the proof of Theorem \ref{thm: height on polynomials} is that while $\mu_f$ can be quite complicated, $\|\mu_f\|_{\lambda_\Ar,1/2}^2$ cannot be too far from the norm of the pull-back of the Arakelov measure $f^*\lambda_\Ar = (f^* \lambda_{\Ar, v})_{v \in M_K}$. It turns out that $\frac{d}{2}\|d^{-1}f^*\lambda_\Ar \|_{\lambda_\Ar,1/2}^2$ is a Weil height on $\Rat_d(K)$, a result of potential interest in its own right (see Section 5 for details). Combining these two results with the observation that $\| \cdot \|_{\lambda,t}^2$ is never too far from $\|\cdot \|^2_{\lambda_\Ar,1/2}$ (Proposition \ref{prop: ag inner product essentially standard inner product}) gives the theorem. \\



Suppose we have two rational maps $f,g: \P^1 \to \P^1$  of degree at least two. The Arakelov--Zhang pairing (also called the dynamical height pairing), denoted here by $(f, g)_{\AZ}$, is a symmetric, non-negative pairing which in some sense is a dynamical distance between $f$ and $g$. It is closely related to the Call--Silverman canonical heights of $f$ and $g$: the pairing $(f, g)_{\AZ}$ vanishes precisely when the canonical heights of $f$ and $g$ are equal, or, equivalently, when $f$ and $g$ have the same canonical adelic measure (see \cite{petsche2012dynamical}, Theorem 3). To some extent, the Arakelov--Zhang pairing  measures how similar the $v$-adic Julia sets of two rational maps $f$ and $g$ are at all places of a number field.\\

The Arakelov--Zhang pairing is not a metric on rational maps; however, Fili (Theorem 9, \cite{fili2017metric}) has discovered that it is expressible as the square of a metric on adelic measures. Specifically, Fili showed that the Arakelov-Zhang pairing is equal to the \textit{mutual energy pairing}: 
\[(f,g)_{\AZ} = \frac{1}{2}(\!(\mu_f - \mu_g, \mu_f - \mu_g)\!),\] 
where
\[ (\!(\mu, \nu)\!) = \sum_{v \in M_K} r_v \iint_{\Aberk\times \Aberk \setminus \Diag_v} \log \frac{1}{\delta_v(x,y)}\; d\mu_{v}\otimes \nu_v(x,y).  \]

Here, $\delta_v(x,y)$ is a continuous extension of $|x-y|_v$ on $\C_v\times \C_v$ to $\Aberk\times \Aberk$ and  $\Diag_v$ is the diagonal of $\C_v \times \C_v$ in $\Aberk \times \Aberk$. (For Archimedean $v$, $\delta_v(x,y)$ is the usual affine distance $|x-y|_v$.) This formulation of the dynamical height pairing has been useful in studying unlikely intersection questions; see, for example, \cite{demarco2020uniform}. The mutual energy pairing itself was introduced by Favre and Rivera-Letelier in connection with a quantitative form of the adelic equidistribution theorem \cite{favre2006equidistribution}. Our inner product could be understood as a mutual energy pairing relative to a fixed adelic probability measure $\lambda$. \\

While quite closely related, our inner product and the mutual energy pairing  are distinct on $\Aad$. For instance, the mutual energy pairing $(\!(\cdot,\cdot)\!)$ is not an inner product on adelic measures of arbitrary mass. Indeed, it follows from a result of Baker--Hsia that if $f$ is a  polynomial defined over $K$, then $(\!(\mu_f, \mu_f)\!) = 0$ (Theorem 4.1 of \cite{baker2005canonical}). However, the mutual energy pairing and our inner product do agree on adelic measures of total mass zero. We may then consider our inner product to be an extension of the mutual energy pairing from the space of adelic measures of total mass zero to $\Aad$. It therefore follows from Theorem 9 of \cite{fili2017metric} and Theorem 1 of \cite{petsche2012dynamical} that if $\mu_f$ and $\mu_g$ are the canonical adelic measures of two rational maps $f$ and $g$ then

\begin{equation}\label{eq: norm equals az pairing}
    \frac{1}{2}\|\mu_f - \mu_g\|_{\lambda,t}^2 = (f, g)_{\AZ}.
\end{equation}

It is a straightforward consequence of Theorem \ref{thm: height on polynomials} and (\ref{eq: norm equals az pairing}) that, when $g$ is fixed, the Arakelov--Zhang pairing is commensurate with the height of $f$. 

\begin{cor}\label{cor: AZ Pairing is a height}
Let $h$ be a Weil height on $\Rat_d(K)$. Fix a rational map $g$ defined over $K$ with degree at least $2$. For any rational map $f$ defined over $K$ with $\deg(f) = d \geq 2$, there are positive constants $c_5, c_6, c_7, c_8$  so that 
\[  c_5 h(f) - c_6 \leq d(f, g)_{\AZ} \leq c_7 h(f) + c_8,   \]
where $c_5, c_7$ depend only on $d$ and tend to $1$ as $d \to \infty$, and where $c_6, c_8$ depend on both $d$,$g$, and $h$. If $f$ is a monic polynomial then 
\[  d(f, g)_{\AZ} = h(f) + O(1)  \]
where the implied constants depend on $g$, $d$, and $h$. 
\end{cor}

Bounds on the Arakelov--Zhang pairing of certain families of rational maps with the power map $z^2$ are studied by Petsche, Szpiro, and Tucker; one such bound, which is generalized by Corollary \ref{cor: AZ Pairing is a height}, is $(z^2, z^2 + c)_{\AZ} =  \frac{1}{2}h(c) +  O(1)$ \cite{petsche2012dynamical}. The lower bound $(z^2, f)_{\AZ} \geq h_f(0)$, where $h_f$ is the Call--Silverman canonical height of a rational map $f$, has been shown by Bridy and Larson \cite{bridy2021arakelov}. Perhaps the strongest known comparisons between the dynamical height pairing and a height on rational maps have been shown by DeMarco, Krieger, and Ye. In \cite{demarco2020uniform}, they show that $(f_{t_1}, f_{t_2})_{\AZ} \asymp h(t_1, t_2)$ for a certain family of Lat\`es maps $f_t$ associated to an elliptic curve written in Legendre form with parameter $t$. In a similar spirit, they show that $(z^2 + c_1, z^2 + c_2)_{\AZ} \asymp h(c_1, c_2)$ \cite{demarco2019common}. These estimates are stronger than those of Corollary \ref{cor: AZ Pairing is a height}, though they hold in specific families of rational maps. It is immediate from Corollary \ref{cor: AZ Pairing is a height} and the Northcott property of heights that $(f,g)_{\AZ}$ can be bounded from below by constant depending on one of the rational maps. 

\begin{cor}\label{cor: lower bound on AZ pairing}
Fix a rational map $f$ of degree $d \geq 2$, defined over $K$. There is a positive constant $C_f = C(f, d, K)$ depending on $f,d$ and $K$ so that for any rational map $g$ defined over $K$ with degree $d$,  $0< (f, g)_{\AZ} \implies C_f \leq  (f, g)_{\AZ}$.
\end{cor}

A uniform lower bound on $(z^2 + c_1, z^2 + c_2)_{\AZ}$ with $c_1 \neq c_2$ and  $c_1, c_2 \in \bar{\Q}$, is given by DeMarco, Krieger, and Ye  (\cite{demarco2019common}, Theorem 1.6). They  use this lower bound on the Arakelov--Zhang pairing to prove a uniform upper bound on the number of preperiodic points shared by $z^2 + c_1$ and $z^2 + c_2$. It is to be hoped that the lower bound of Corollary \ref{cor: lower bound on AZ pairing} has analogous implications for the number of common preperiodic of two rational maps of degree $d$, but we do not explore this possibility here. 

\begin{rem}\label{rem: small points intro}
The estimates of Theorems \ref{thm: norm of monic poly is a height} and \ref{thm: height on polynomials} give estimates on $\|\mu\|_{\lambda, t}^2$ when $\mu = \mu_f$ is the canonical adelic measure associated to a rational map $f$. It is natural to wonder if the norm of a general adelic measure can be estimated, or at least bounded from below. It is a formal consequence of the definitions that $\|\mu\|_{\lambda,t}^2 \geq t$ for any adelic probability measure $\mu$; see Corollary \ref{cor: pos def of ag inner product}. There is a class of adelic probability measures for which better lower bounds are possible. We say that an adelic measure $\mu$ has \textit{points of small height} if there is a sequence of distinct point $\{\alpha_n\}_{n=1}^\infty$ for which $h_\mu(\alpha_n) \to 0$ as $n \to \infty$, where $h_\mu$ is the Favre--Rivera-Letelier height of $\mu$. For example, every canonical adelic measure $\mu_f$ associated to a rational map of degree at least $2$ has small points (indeed, $h_{\mu_f}(\alpha) =0$ if and only if $\alpha$ is preperiodic).  For such measures, it is straightforward to show that $\|\mu\|^2_{\lambda_{\Ar, 1/2} }\geq \frac{1}{2}\log(2)$ with equality if and only if $\mu = \mu_{z^2}$. See Proposition \ref{prop: small points} for details. 
\end{rem}

\subsection{Outline}
We discuss the notation and tools needed in the next section. Section \ref{sec: adelic measures} contains the proof of Theorem \ref{thm: inner product into}. Section \ref{sec: comparison w ar height MONIC} has the proof of Theorem \ref{thm: heigh on monic poly}. Section \ref{sec: comparison w ar ht rat} contains the proof of \ref{thm: height on polynomials} and of Corollary \ref{cor: lower bound on AZ pairing}. Section \ref{sec: examples} contains some examples where the norm can be computed or bounded explicitly. Here we obtain exact expressions for the norm (with respect to the Arakelov measure) associated to Chebyshev polynomials, and for the norm associated to power maps after an affine change of variables. We also see that the norm associated to a  certain family of Latt\`es maps can be bounded in a stronger way than the general bounds of Theorem \ref{thm: height on polynomials} for rational maps. 

\subsection{Acknowledgements} Many of the results of this article have appeared in my PhD thesis \cite{oberlythesis}, which was written at Oregon State University. To my advisor, Clayton Petsche, I owe my deepest thanks. I am also very grateful to the anonymous referee for several insightful comments and useful suggestions. 
\section{Preliminaries}\label{sec: preliminaries}

Fix a number field $K$ and let $M_K$ be the set of places of $K$. Let $\C_v$ be the completion of a fixed algebraic closure of $K_v$. We denote by $|\cdot|_v$ the absolute value on $\C_v$ whose restriction to $\Q$ coincides with the usual real or $p$-adic absolute value. Let $r_v = [K_v: \Q_v]/[ K :\Q]$ be the ratio of local-to-global degrees. It it well known that $K$ satisfies the \textit{product formula}, which states that $\prod_{v \in M_K} |x|_v^{r_v} = 1$ for all $x \in K^\times$. Moreover the \textit{extension formula} reads $\sum_{\substack{w \in M_L\\w \mid v}} r_w = r_v$, where $L$ is a finite extension of $K$, and where $w \mid v$ means that $w$ lies above $v$ (so $|\alpha|_w = |\alpha|_v$ for all $\alpha \in K$). We write $M_K^0$ for the set of non-Archimedean places of $K$ and $M_K^\infty$ for the Archimedean places.  \\

Of central importance to this article is the \textit{v-adic chordal metric}, $||\cdot,\cdot ||_v: \P^1(\C_v) \times \P^1(\C_v) \to [0, 1]$. In homogeneous coordinates $x = [x_1 : x_2]$ and $y = [y_1 : y_2]$ the chordal metric is given by
\[ ||x,y||_v =  \begin{dcases}
                       \frac{ |x_1 y_2 - x_2 y_1|_v }{\sqrt{ |x_1|_v^2 + |x_2|_v^2}\sqrt{ |y_1|_v^2+ |y_2|_v^2}}, & \qif v\mid \infty;\\
                       \frac{|x_1 y_2 - x_2y_1|_v}{\max(|x_1|_v, |x_2|_v) \max( |y_1|_v, |y_2|_v)} & \qif v \nmid \infty.
                       \end{dcases}
\]
It is well known that $||\cdot, \cdot||_v$ is a metric on $\P^1(\C_v)$, and that when $v$ is not Archimedean, then $||\cdot, \cdot||_v$ satisfies the strong triangle inequality (see, for example, chapters 1 and 3 of \cite{silverman2007arithmetic}). For non-Archimedean $v$, there is a continuous extension of $||\cdot, \cdot||_v$ to $\Pberk$ called the \textit{spherical kernel} which we denote again by $||\cdot, \cdot||_v$. The spherical kernel is not a metric on $\Pberk$; indeed, the Gauss point $\zeta_{\Gauss}$ satisfies $||z, \zeta_\Gauss||_v = 1$ for all $z \in \Pberk$. See Proposition 4.7 of  \cite{baker2010potential} for more properties of the spherical kernel.\\


There is an analogous extension of the absolute value $|x-y|_v$ on $\C_v$ to the \textit{Berkovich affine line}, $\Aberk$, which, following \cite{baker2010potential}, we refer to as the \textit{Hsia Kernel}, which we will write as $\delta_v(x,y)$. (In \cite{baker2010potential}, this is written $\delta(x,y)_\infty$ to distinguish the Hsia kernel from the generalized Hsia kernel. In the work of Favre and Rivera-Letelier, this extension is written $\sup\{ S, S'\}$.) 
Similar to the classical relationship between  the affine distance on $\C$ and the chordal metric, we have 
\begin{equation}\label{eq: chordal to affine dist}
    ||x,y||_v =  \delta_v(x,y) ||x,\infty||_v  ||y, \infty||_v
\end{equation}
for any two points $x, y\in \Aberk \cong \Pberk \setminus \{\infty\}$. We will also write $\delta_v(x,y) = |x-y|_v$ for the affine distance on $\C_v$ when $v$ is Archimedean.  

\subsection{The Measure-valued Laplacian}

Given a sufficiently regular function $f: \Pberk \to \R$, the \textit{measure-valued Laplacian} is a an operator $\Delta$ which assigns to $f$ a Radon measure $\Delta f$ on $\Pberk$ of total mass $0$ (Chapter 5, \cite{baker2010potential}). When $v$ is Archimedean, we identify $\Pberk$ with the Riemann Sphere $\P^1(\C) \cong \C \cup \{\infty\}$. In this case, the measure valued Laplacian is defined for twice continuously differentiable functions $f: \P^1(\C) \to \R$ by $\Delta f(z) = \frac{-1}{2\pi}(\partial_{xx} + \partial_{yy})f(z)\; dx\; dy$,
where $z = [z: 1]$ and $z = x + iy$. We then extend $\Delta$ in a distributional sense.  When $v$ is non-Archimedean, the measure-valued Laplacian is first defined on finitely branching subgraphs of $\Pberk$, and then extended to functions $f: \Pberk \to \R$ which satisfy appropriate coherence conditions. Following \cite{baker2010potential}, we denote by $\BDV_v$ the space of functions of \textit{bounded differential variation} on $\Pberk$; in a certain sense, this is the largest class of functions for which the measure-valued Laplacian can be defined. For Archimedean $v$, we will also use the notation $\BDV_v$ to denote those functions whose distributional Laplacian also defines a Radon measure; in particular, this class contains the sub-harmonic functions (\cite{ransford1995potential}, Theorem 3.7.2). We record the standard properties of the Laplacian in the following proposition.

\begin{prop}\label{prop: properties of the laplacian}
Fix $v \in M_K$ and let $f, g: \Pberk \to \R$ be in $C(\Pberk) \cap \BDV_v$. The measure valued Laplacian has the following properties:
\begin{enumerate}
    \item The Laplacian is self-adjoint: $\int f\; d\Delta g = \int g\; d \Delta f$ whenever $f$ is $\Delta g$-integrable and $g$ is $\Delta f$-integrable; 
    \item if $\Delta f= \Delta g$, then $f$ and $g$ differ by a constant;
    \item $\int f\; d\Delta f \geq 0$ with equality if and only if $f$ is constant.
\end{enumerate}
\end{prop}

When $v$ is Archimedean, these are standard facts which follow from Green's identities. See Corollary 5.38 of \cite{baker2010potential} for a proof of the non-Archimedean case. 


\subsection{Local Potential Theory} 


Given a complex Radon measure $\mu$ on $\Pberk$, the \textit{(local) potential of} $\mu$ is the function $U_\mu:\Pberk \to \C$ defined by
\[ U_\mu(x) = \int_{\Pberk} \log\frac{1}{||x, y||_v} \; d\mu(y),   \]
 supposing that the integrand is $\mu$-integrable for all $x$. We say that $\mu$ \textit{has continuous potential} if $U_\mu$ is continuous. For each place $v \in M_K$ set 

\begin{equation}\label{eq: definition of arakelov measure}
d\lambda_{\Ar, v} = \begin{cases}
\frac{d\ell(z)}{\pi(1 + |z|_v^2)^2} &\qif v \mid \infty;\\
 d\delta_{\zeta_\Gauss}            &\qif v \nmid \infty,
\end{cases}
\end{equation}
where $\ell$ is the Lebesgue measure on $\C$ and $\delta_{\zeta_\Gauss}$ is the point mass concentrated at the Gauss point $\zeta_\Gauss$. For Archimedean $v$, a straightforward computation shows that $\lambda_{\Ar, v}$ is invariant under the action of the unitary group $U_2(\C)$ on $\P^1(\C)$. 

\begin{prop}\label{prop: laplacian and spherical potentials}
Let $v \in M_K$ and let $\mu$ be a signed Radon measure on $\Pberk$ with continuous potential $U_\mu$. Then $U_{\mu} \in \BDV_v$ and $\Delta U_\mu = \mu - \mu(\Pberk) \lambda_{\Ar, v}$.
\end{prop}

\begin{proof}
For non-Archimedean $v$, the proposition is Examples 5.19 and 5.22 of \cite{baker2010potential}. For Archimedean $v$, Theorem 1.5 of \cite{lang2012introduction} states that $\Delta_x \log||x, y||_v^{-1} = \delta_y - \lambda_{\Ar, v}$. Thus, if $\phi$ is a smooth test function, then
\begin{align*}
    \int U_\mu(x)\, d\Delta\phi(x) &= \int \int \log \frac{1}{||x, y||_v}\, d\mu(y)\, d\Delta\phi(x)\\
                                   &= \int \int \log \frac{1}{||x, y||_v}\, d\Delta \phi(x)\, d\mu(y)\\
                                   &= \int \left( \phi(y) - \int \phi(x) \, d\lambda_{\Ar, v}(x) \right) \, d\mu(y)\\
                                   &= \int \phi(y)\, d\mu(y) - \mu(\Pberk)\int \phi(x)\, d\lambda_{\Ar, v}(x),
\end{align*}    
where the use of Fubini's theorem in the second line is justified as $U_\mu$ is continuous, hence $|\Delta \phi|$-integrable, and where the distributional equation $\Delta_x \log ||x,y||_v^{-1} = \delta_y(x)-\lambda_{\Ar, v}(x)$ is used in the third line. Consequently $\Delta U_\mu = \mu - \mu(\Pberk)\lambda_{\Ar, v}$ for Archimedean $v$ as well. 
\end{proof}

Given a probability measure on $\Pberk$ with continuous potential $U_\mu$, an \textit{Arakelov--Green's} function of $\mu$ is a map $g_\mu :\Pberk \times \Pberk \to (-\infty, \infty]$ which is symmetric, continuous off of the diagonal of $\Pberk \times \Pberk$ and satisfies the distributional equation $\Delta_x g_\mu(x, y) = \delta_y - \mu$. These conditions determine $g_{\mu}$ up to an additive constant. In terms of the potentials $U_{\mu}$, $g_{\mu}$ is given by the expression
\begin{equation}\label{eq: greens function in chordal}
    g_\mu(x, y) = \log \frac{1}{||x, y||_v} - U_\mu(x) - U_\mu(y) + C 
\end{equation}  
where $C$ is a constant (p. 241 of \cite{baker2010potential}). For example, an Arakelov--Green's function of $\lambda_{\Ar, v}$ is $g_{\lambda_{\Ar}, v}(x,y) = \log ||x, y||_{v}^{-1} + C$, a fact which follows from the distributional equation $\Delta_x \log ||x, y||_v^{-1} = \delta_y - \lambda_{\Ar, v}$. 
\section{Adelic Measures}\label{sec: adelic measures}

We come to the main object of study in this article. We adopt the convention that complex or real valued measures have finite total variation.

\begin{defn}
An \textit{adelic measure} defined over $K$ is a sequence $\mu = ( \mu_v)_{v \in M_K}$ indexed by places $v \in M_K$ so that:
\begin{enumerate}
    \item For each $v\in M_K$, $\mu_v$ is a complex-valued Radon measure on $\Pberk$;
    \item There is a constant $c_\mu \in \C$ so that $\mu_v(\Pberk) = c_\mu$ for all $v \in M_K$;
    \item For all but finitely many $v$, $\mu_v = c_\mu \delta_{\zeta_\Gauss}$;
    \item  $\mu_v$ has continuous potentials for all $v \in M_K$.
\end{enumerate}
We say that an adelic measure $\mu$ is a signed, probability or positive adelic measure if each $\mu_v$ is a signed, probability or positive measure respectively.  We write $U_{\mu, v}$ instead of the more cumbersome $U_{\mu_v}$ for the local potential of $\mu_v$. Likewise we will write $g_{\mu,v}$ for the Arakelov--Green's function $g_{\mu_v}$ of $\mu_v$. The \textit{total mass} of an adelic measure is defined to be the constant $c_\mu$. We denote by $\Aad$ the complex vector space of all proper complex adelic measures with continuous potentials, and by $\Aad^0$ the subspace consisting of those adelic measures with total mass $0$.  
\end{defn}

Observe that if $\mu = (\mu_v)_{v \in M_K}$ is an adelic measure with $\mu_v = \mu_{r, v} + i\mu_{i, v}$ for signed measures $\mu_{r, v}, \mu_{i, v}$ at each place, then $(\mu_{r, v})_{v \in M_K}$ and $(\mu_{i,v})_{v \in M_K}$ are again (signed) adelic measures. Indeed, $\Re U_{\mu, v} = U_{\mu_{r}, v}$ so that $\mu_{r, v}$ has continuous potentials at each place $v$; moreover $\Re c_\mu = \Re \mu_v(\Pberk) = \mu_{r, v}(\Pberk)$ so that $\mu_{r, v}(\Pberk)$ is a constant independent of $v$. So $\mu_r = (\mu_{r, v})$ is a signed adelic measure. Similar remarks apply to $(\mu_{i, v})_{v \in M_K}$. \\

The primary examples of adelic probability measures in $\Aad$ relevant to arithmetic dynamics are the adelic canonical measures associated to a rational map, which can be defined as the weak (really, weak-$*$) limit of the sequence $(d^{-1} f^*)^n \lambda_{\Ar,v}$ of normalized pull-backs of the Arakalov measure (or, any other measure with continuous potentials). It follows from this definition that $f^* \mu_{f,v} = d\mu_{f,v}$; in fact, $\mu_{f,v}$ is the unique log-continuous  probability measure satisfying this equation (Theorem (d) in \cite{freire1983invariant}, Th\'eor\`eme A in \cite{favre2010theorie}).  As shown in \cite{favre2004theoreme}, $\mu_{f,v} = \delta_{\zeta_\Gauss}$ precisely when $f$ has good reduction, and so $\mu_{f, v} = \delta_{\zeta_\Gauss}$ for all but finitely many $v$. That $\mu_{f}$ has continuous potentials (in fact, Holder continuous) is also due to Favre and Rivera-Letelier for non-Archimedean $v$, and due to Ma{\~n}{\'e} for Archimedean $v$ \cite{mane2006hausdorff}. So $\mu_f = (\mu_{f, v})_{v \in M_K}$ is indeed an  element of $\Aad$. \\

\subsection{A Family of Inner Products} Fix an adelic probability measure $\lambda = (\lambda_v)_{v \in M_K} \in \Aad$. At each place $v$, let $g_{\lambda, v}(x, y)$ be the Arakelov--Green's function of $\lambda_v$, normalized so that 
\begin{equation}\label{eq: normalization of green's function}
\int g_{\lambda, v}(x, y)\, d\lambda(x) =: t_v \geq 0
\end{equation}
at each place, and where $t_v = 0$ for all but finitely many $v$, and $\sum_{v \in M_K} r_v t_v =: t > 0$. Given adelic measures $\mu = (\mu_v)_{v \in M_K}, \nu = (\nu_v)_{v \in M_K} \in \Aad$, define 
\[ \< \mu_v, \nu_v\>_{\lambda, t, v} = \int\int g_{\lambda, v} (x,y)  \, d\mu_v(x)\, d\overline{\nu_v}(y)\] 
and
\begin{equation}\label{eq: defn of ag inner product}
    \< \mu, \nu\>_{\lambda, t} = \sum_{v \in M_K} r_v \< \mu_v, \nu_v\>_{\lambda, t, v}. 
\end{equation}  

Note that $\iint g_{\lambda, v} \, d\mu_v d\overline{\nu_v}$ exists as $\mu_v, \nu_v$ have continuous potentials. Moreover, that $g_{\lambda, v}(x, y) = \log ||x, y||_v^{-1}$ and $\mu_v = c_\mu \delta_{\zeta_\Gauss}$ for all but finitely many $v$, and also that $\int \log||x, y||_v^{-1} \, d\delta_{\zeta_{\Gauss}}(y) = 0$ implies $\< \mu_v, \nu_v\>_{\lambda, t, v} = 0$ for all but finitely many places $v$. So the sum (\ref{eq: defn of ag inner product} defining $\< \mu, \nu\>_{\lambda, t}$ is finite.  Notice that $\<\cdot,\cdot\>_{\lambda,t,v}$ is linear in the first coordinate and is conjugate-symmetric; consequently so is $\<\cdot,\cdot \>_{\lambda, t}$. \\

We write $\|\mu_v\|_{\lambda,  t,v}^2 = \< \mu_v, \mu_v\>_{\lambda, t, v}$ and $\|\mu\|^2_{\lambda, t} = \< \mu, \mu\>_{\lambda, t}$. It follows from our choice of normalization of $g_{\lambda, v}$ that
\[ 
\| \lambda_v\|^2_{\lambda, t,v} = \int\int g_{\lambda, v}(x,y)\, d\lambda_v(x)\, d\lambda_v(y) = t_v
\]
and so
\[\|\lambda\|_{\lambda, t}^2 = \sum_{v \in M_K} r_v t_v  = t. \]

Our first theorem describes the properties of this pairing and shows that it is an extension of the mutual energy pairing investigated by Favre and Rivera-Letelier, and by Fili. We adopt the terminology common in probability and say that a sequence of (signed or complex-valued) Radon measures $\{\mu_n\}_{n=1}^\infty$ converges weakly to a measure $\mu$ if $\int \phi\, d\mu_n \to \int \phi \,d\mu$ for all continuous $\phi$. 

\begin{thm}\label{thm: AG inner product}
    Let $\lambda = (\lambda_v)_{v \in M_K} \in \Aad$ be an adelic probability measure, and, at each place $v$, choose a normalization of $g_{\lambda, v}$ as in (\ref{eq: normalization of green's function}). Then $\< \cdot, \cdot\>_{\lambda, t}$ defines an inner product on $\Aad$ which agrees with the mutual energy pairing on $\Aad^0$. Moreover, if $\mu_n, \mu \in \Aad$ with $\| \mu_n - \mu\|_{\lambda, t} \to 0$ and if $\sup_{n \geq 1} |\mu_{n,v}|(\Pberk) < \infty$, then $\mu _{n, v}$ converges weakly to $\mu_v$.
\end{thm}

\subsection{Proof of  Theorem \ref{thm: AG inner product}}Rather than studying the whole family of pairings $\< \cdot, \cdot\>_{\lambda, t}$, it will be convenient to study one particular pairing. Proposition \ref{prop: ag inner product essentially standard inner product} shows that there is little loss of generality in doing so. 

Given $\mu = (\mu_v)_{v \in M_K}, \nu = (\nu_v)_{v \in M_K}$ in $\Aad$, set
\begin{equation*}
     \< \mu_v, \nu_v\>_v = \int_{\Pberk}\int_{\Pberk} \log\frac{1}{||x, y||_v} \;d\mu_v(x) d\overline{\nu_v}(y) = \int_{\Pberk} U_{\mu,v}(y)\; d\overline{\nu_v}(y),
\end{equation*}
and 
\begin{equation}\label{eq: innerproddef} 
\< \mu, \nu\> = \sum_{v \in M_K} r_v \< \mu_v, \nu_v\>_v. 
\end{equation}
We will write $\|\mu_v\|^2_v = \< \mu_v, \mu_v\>_v$ and $\|\mu\|^2= \< \mu, \mu\>$. 

This is the pairing associated to $\lambda_\Ar = (\lambda_{\Ar, v})_{v \in M_K}$ with $t = 1/2$. 
Thus
\[ \< \mu_v, \nu_v\>_v = \< \mu_v, \nu_v\>_{\lambda_\Ar, 1/2, v}. \]
To see this, note that $g_{\lambda_{\Ar}, v}(x,y) = \log ||x, y||_v^{-1} + C$. For non-Archimedean $v$, $\lambda_{\Ar, v} = \delta_{\zeta_{\Gauss}}$ and $\int \log ||x, y||_v^{-1} \, d\delta_{\zeta_\Gauss}(x) = \log ||\zeta_\Gauss, y||_v^{-1} = 0$. For Archimedean $v$, a straightforward computation, using the invariance of both the projective metric $||\cdot, \cdot||_v$ and $\lambda_{\Ar, v}$ under the action of $U_2(\C)$ on $\P^1(\C)$, shows that $\int \log ||x, y||_v^{-1} \, d\lambda_{\Ar, v}(x) = \frac{1}{2}$. Therefore $\sum_{v \in M_K}r_v \int\int g_{\lambda_\Ar, v} \, d\lambda_{\Ar, v}d \lambda_{\Ar, v} = \sum_{v \mid \infty} r_v \frac{1}{2} =\frac{1}{2}$, where we have used the fact that $\sum_{v \mid \infty} r_v = 1$.

\begin{prop}\label{prop: ag inner product essentially standard inner product}
Let $\lambda = (\lambda_v)_{v \in M_K}$ be an adelic probability measure, and choose a normalization of Arakelov--Green's functions $g_{\lambda,v}$ as in (\ref{eq: normalization of green's function}). Then 
\begin{enumerate}
    \item $\<\mu, \nu\>_{\lambda, t} = \< \mu - c_\mu \lambda, \nu - c_\nu \lambda\> + c_\mu \overline{c_\nu} t$  for any $\mu, \nu \in \Aad$;\\
    \item  there is a constant $\kappa$ depending only on $\lambda$ so that $|\< \mu, \nu\>_{\lambda, t} - \< \mu, \nu\>| \leq c_\mu c_\nu(\kappa+t)$ whenever $\mu$ and $\nu$ are positive adelic measures.
\end{enumerate}
\end{prop}

We remark that item (1) shows that $\< \cdot, \cdot\>_{\lambda,t}$ depends only on $\sum_{v \in M_K} r_v t_v = t$ and not on the specific normalization of $g_{\lambda, v}$ at each place. 

\begin{proof}
It is a formal consequence of (\ref{eq: greens function in chordal}) that
\begin{align*}
    t_v &=  \int \int g_{\lambda, v}(x,y) \, d\lambda_v(x) \, d\lambda_v(y) \\
        &= \int \int \left\{ \log \frac{1}{||x, y||_v} - U_{\lambda,v}(x) - U_{\lambda,v}(y) + C \right\} \, d\lambda_v(x) d\, \lambda_v(y)\\
        &= -\|\lambda_v\|_v^2 + C,
\end{align*} 
and so 
\begin{equation}\label{eq: non normalized green's functions}
    g_{\lambda, v}(x,y) = \log\frac{1}{||x, y||_v} - U_{\lambda, v}(x) - U_{\lambda, v}(y) + \|\lambda_v\|_v^2 + t_v. 
\end{equation}
Note that $\int U_{\lambda,v}(x) \, d\mu_v(x) = \int \int \log ||x, y||_v^{-1} \, d\lambda_v(y)\, d\mu_v(x) = \< \mu_v, \lambda_v\>_v$, and similarly that $\int U_{\lambda, v}(y)\, \overline{\nu_v}(y) = \< \lambda_v, \nu_v\>_v$. By (\ref{eq: non normalized green's functions}), we have
    \begin{align*}
        \< \mu_v, \nu_v\>_{\lambda, t, v} &= \int \int g_{\lambda, v}(x, y)\; d\mu_v(x)\; d\nubar_v(y)\\
                                       &= \int\int \left\{\log \frac{1}{||x, y||_v} - U_{\lambda, v}(x) - U_{\lambda, v}(y) + \|\lambda_v\|^2_v + t_v\right\} \; d\mu_v(x)\; d\nubar_v(y)\\
                                       &= \< \mu_v, \nu_v\>_v - \overline{c_\nu} \< \mu_v, \lambda_v\>_v - c_\mu \< \lambda_v, \nu_v\>_v + c_\mu \overline{c_\nu} \|\lambda_v\|_v^2 + c_\mu \overline{c_\nu} t_v\\
                                       &= \< \mu_v, \nu_v - c_\nu \lambda_v\>_v -  \< c_\mu \lambda_v, \nu_v\>_v + \< c_\mu \lambda_v, c_\nu \lambda_v\>_v  + c_\mu \overline{c_\nu} t_v\\
                                       &= \< \mu_v, \nu_v - c_\nu \lambda_v\>_v -  \< c_\mu \lambda_v, \nu_v - c_\nu \lambda_v \>_v +   c_\mu \overline{c_\nu} t_v\\
                                       &= \< \mu_v - c_\mu \lambda_v, \nu_v - c_\nu \lambda_v\>_v + c_\mu \overline{c_\nu}t_v,
    \end{align*}
    where we have used the linearity and conjugate-symmetry of $\<\cdot, \cdot\>_v$ several times. Multiplying by $r_v$ and summing over all places of $K$ gives item (1). \\
    For (2), note that as $U_{\lambda, v}(x)$ is continuous on $\Pberk$, (\ref{eq: greens function in chordal}) implies that there is a constant $\kappa_v$ so that 
    \begin{equation}\label{eq: diff between green and standard}
      \big|  g_{\lambda, v}(x, y) - t_v  - \log ||x, y||_v^{-1} \big|\leq  \kappa_v
    \end{equation} 
    for all $x \neq y$. Moreover, when $\lambda_{v}=\delta_{\zeta_\Gauss}$, $g_{\lambda,v}(x,y) -t_v= \log||x,y|_v^{-1}$. Therefore $\kappa_v = 0$ for all but finitely many $v\in M_K$. So, if $\mu = (\mu_v)_{v \in M_K}, \nu = (\nu_v)_{v \in M_K}$ are positive adelic measures, then integrating (\ref{eq: diff between green and standard}) at each place gives  \( \big|\< \mu_v, \nu_v\>_{\lambda, t, v} - c_\mu c_\nu t_v - \< \mu_v, \nu_v\> \big| \leq \kappa_v c_\mu c_\nu \). Therefore, using the fact that $t_v$ is non-negative, we have $|\<\mu_v, \nu_v\>_{\lambda, t, v} - \< \mu_v, \nu_v\>|\leq c_\mu c_\nu (\kappa_v + t_v)$. Multiplying by $r_v$ and summing over all places of $K$ gives item (2) with $\kappa = \sum_v r_v \kappa_v$. 
\end{proof}

For the remainder of this section, fix an adelic probability measure $\lambda = (\lambda_v)_{v \in M_K}$ and a normalization of $g_{\lambda,v}$ at each place $v$. The next proposition will allow us to use the positive definiteness of the mutual energy pairing on $\Aad^0$ to deduce the positive definiteness of $\< \cdot, \cdot\>_{\lambda, t}$.

\begin{prop}\label{prop: ag inner product and mep}
For any $\mu, \nu \in \Aad$, we have $\< \mu, \nu\>_{\lambda, t} = (\!( \mu - c_\mu \lambda, \nu - c_\nu \lambda)\!) + c_\mu \overline{c_\nu}t$. In particular, $\< \mu, \nu\>_{\lambda, t} = (\!( \mu, \nu)\!)$ whenever $\mu, \nu \in \Aad^0$. 
\end{prop}

\begin{proof}
Let $\mu = (\mu_v)_{v\in M_K}$ and $\nu = (\nu_v)_{v \in M_K}$ be adelic measures. We assume first that $\mu_v(\Pberk) = \nu(\Pberk) = 0$. That $\mu_v$ and $\nu_v$ are log-continuous Radon measures on $\Pberk$ implies that $\mu$ and $\nu$ do not charge the point at infinity; moreover $\mu_v \otimes \overline{\nu_v}(\Diag) = 0$. By (\ref{eq: chordal to affine dist}), we have

\begin{align*}
    \< \mu_v, \nu_v\>_v &= \int_{\Pberk}\int_{\Pberk}  -\log||x, y||_v\; d\mu_v(x) d\overline{\nu_v}(y)\\
                &= \int_{\Aberk}\int_{\Aberk} \big( -\log \delta_v(x,y) - \log ||x, \infty||_v - \log||y, \infty||_v\big)\; d\mu_v(x)\; d\overline{\nu_v}(y)\\
                &= \iint_{\Aberk\times \Aberk} -\log \delta_v(x,y)\; d\mu_v\otimes\overline{\nu_v}(x,y) - \overline{\nu_v}(\Pberk)\int_{\Pberk} \log ||x, \infty||_v \; d\mu_v(x)\\
                &\qquad - \mu_v(\Pberk) \int_{\Pberk} \log ||y, \infty||_v \; d\overline{\nu_v}(y) \\
                &=  \iint_{\Aberk\times \Aberk \setminus \Diag}  -\log \delta_v(x,y)   \;d\mu\otimes \overline{\nu_v}(x,y)\\
                &= (\!( \mu_v, \nu_v)\!)_v,
\end{align*}
where we have used the assumption $\mu_v(\Pberk) = \nu_v(\Pberk) = 0$. Multiplying by $r_v$ and summing over all places shows that $\< \mu, \nu\> = (\!(\mu, \nu)\!)$ when $\mu$ and $\nu$ have total mass $0$. \\
Now, assuming that $\mu$ and $\nu$ have arbitrary mass $c_\mu$ and $c_\nu$, we apply the above result to $\mu - c_\mu \lambda$ and $\nu - c_\nu \lambda$ and invoke item (1) of Proposition (\ref{prop: ag inner product essentially standard inner product}) to find 
\[ \< \mu, \nu\>_{\lambda, t} - c_\mu \overline{c_\nu} t = \< \mu- c_\mu \lambda, \nu - c_\nu \lambda\> = (\!( \mu - c_\mu \lambda, \nu - c_\nu \lambda)\!).  \]
\end{proof}

\begin{rem}
    If $\lambda$ is taken to be the canonical adelic measure $\mu_g$ of a rational map $g$ with degree at least $2$, and if $f$ is another rational map of degree at least $2$, then Proposition \ref{prop: ag inner product and mep} gives the following relationship between $\|\cdot\|_{\mu_g, t}$ and the Arakelov--Zhang pairing:
\begin{equation}\label{eq: ag norm az pairing}
    \lim_{t\to 0} \| \mu_f\|_{\mu_g, t}^2 = (\!( \mu_f -\mu_g, \mu_f - \mu_g)\!) = 2(f, g)_{\AZ},   
\end{equation}
where, as remarked in the introduction, the last equality is due to Fili (Theorem 9, \cite{fili2017metric}). 
\end{rem}

\begin{cor}\label{cor: pos def of ag inner product}
    For all $\mu\in \Aad$, $\< \mu, \mu\>_{\lambda, t} = \|\mu\|^2_{\lambda, t} \geq |c_\mu|^2 t$ with equality if and only if $\mu = c_\mu \lambda$. 
\end{cor}

\begin{proof}
Propositions 2.6 and 4.5 of \cite{favre2006equidistribution} state that if $\rho_v$ is a signed measure of total mass $0$ on $\Pberk$, then $(\!( \rho_v, \rho_v)\!)_v \geq 0$, with equality if and only if $\rho_v = 0$. See also Theorem 8.72 of \cite{baker2010potential} and Theorem 5.3 of \cite{baker2006equidistribution}. Consequently, if $\mu = (\mu_v)_{v \in M_K} \in \Aad$ is a signed adelic measure then Proposition \ref{prop: ag inner product and mep} implies \(\< \mu, \mu\>_{\lambda, t} = (\!(\mu - c_\mu \lambda, \mu - c_\mu \lambda)\!) + |c_\mu|^2 t \geq  |c_\mu|^2 t \) with equality if and only if $\mu = c_\mu \lambda$. 
We recall that if $\mu$ is a complex-valued adelic measure, say $\mu_v = \mu_{r, v} + i\mu_{i, v}$ at each place, then $\mu_r = (\mu_{r, v})_{v \in M_K}$ and $\mu_i = (\mu_{i, v})_{v \in M_K}$ are signed adelic measures. Therefore
\begin{align*}
    \< \mu_v, \mu_v\>_{\lambda, t} &=\< \mu_{r} + i\mu_{i}, \mu_{r} + i\mu_{i} \>_{\lambda, t} \\
    &= \|\mu_{r}\|_{\lambda, t}^2 - \< \mu_{r}, i\mu_{i}\>_{\lambda, t} - \< i\mu_{i}, \mu_{r}\>_{\lambda, t} + \|\mu_{i}\|^2_{\lambda, t} \\
    &= \|\mu_{r}\|_{\lambda, t}^2  + i \< \mu_{r}, \mu_{i}\>_{\lambda, t} - i\<\mu_{i}, \mu_{r}\>_{\lambda, t} + \|\mu_{i}\|^2_{\lambda, t}\\
    &= \|\mu_{r}\|^2_{\lambda, t} + \|\mu_{i}\|^2_{\lambda, t}.
\end{align*} 
Now $\|\mu_r\|_{\lambda, t}^2 \geq |c_{\mu_r}|^2 t$ with equality if and only if $\mu_r = c_{\mu_r} \lambda$; similar remarks apply to $\mu_i$. Thus
\[ \|\mu\|^2_{\lambda, t} = \|\mu_r\|^2_{\lambda, t} + \|\mu_i\|^2_{\lambda, t} \geq (|c_{\mu_r}|^2 +|c_{\mu_i}|^2 )t = |c_\mu|^2 t \]
with equality if and only if $\mu_r = c_{\mu_r}\lambda$ and $\mu_i = c_{\mu_i}\lambda$. That is, equality holds if and only if $\mu = c_\mu \lambda = (c_{\mu_r} + ic_{\mu_i})\lambda$.
\end{proof}

It follows from Corollary \ref{cor: pos def of ag inner product} that $\|\mu\|^2_{\lambda, t} \geq 0$ with equality if and only if $\mu = 0\lambda = 0$. Thus $\<\cdot, \cdot\>_{\lambda, t}$ is positive-definite on $\Aad$. Since $\< \cdot, \cdot \>_{\lambda, t}$ is linear in the first coordinate and conjugate symmetric, then we have shown that $\< \cdot, \cdot\>
_{\lambda, t}$ is an inner product on $\Aad$. Proposition \ref{prop: ag inner product and mep} shows that $\< \cdot, \cdot\>_{\lambda, t}$ extends the mutual energy pairing to $\Aad$. To finish the proof of Theorem \ref{thm: AG inner product} we need only show the statement on weak convergence. We remark that Favre and Rivera-Letelier show a more general statement which holds for differences of probability measures; see Propositions 2.11 and 4.12 of \cite{favre2006equidistribution}. We denote by $|\mu|(\Pberk)$ the total variation of a complex valued measure $\mu$ on $\Pberk$. 

\begin{prop}
Let $\mu_n, \mu \in \Aad$ and suppose that $\|\mu_n - \mu\|^2_{\lambda,t} \to 0$ as $n \to \infty$. Let $v \in M_K$. Then: 
\begin{enumerate}
    \item For every $\phi \in C(\Pberk) \cap \BDV_v$, $\int \phi\, d\mu_{n, v} \to \int \phi\, d\mu_v$ as $n \to \infty$;
    \item If $\sup_{n \geq 1} |\mu_{n, v}|(\Pberk) < \infty$ then $\mu_{n,v} \to \mu_{v}$ weakly as $n \to \infty$.
\end{enumerate}
\end{prop}

\begin{proof}
By Replacing $\mu_n$ with $\mu_n- \mu$, we assume that $\mu = 0$. We note that it suffices to prove the proposition with $\|\cdot\|^2$ in place of $\|\cdot\|_{\lambda,t}^2$. Indeed, by Proposition \ref{prop: ag inner product essentially standard inner product}, if $\|\mu_n\|_{\lambda, t}^2 \to 0$ then $\|\mu_n - c_{\mu_n}\lambda\|^2 + |c_{\mu_n}|^2t \to 0$ as $n\to \infty$; so $\|\mu_n - c_{\mu_n}\lambda\|$ and $|c_{\mu_n}|$ both tend to $0$ as $n$ goes to infinity. Then the triangle inequality implies  $\limsup_{n \to \infty} \|\mu_n\| \leq \limsup_{n\to \infty} (\|\mu_n - c_{\mu_n}\lambda\| + |c_\mu| \|\lambda\|) = 0.$ So $\|\mu_n\|_{\lambda,t} \to 0$ implies $\|\mu\| \to 0$. 

Let $E_v = C(\Pberk) \cap \BDV(\Pberk)$, and let $\phi \in E_v$. In both non-Archimedean and Archimedean cases, $E_v$ is a dense subspace of $C(\Pberk)$ (Proposition 5.4 and Example 5.18 of \cite{baker2010potential} when $v \nmid \infty$; for $\P^1(\C)$, this is well known - in fact, $E_v$ includes the smooth functions for Archimedean $v$). \\
That $\phi \in E_v$ implies that $\Delta \phi$ exits, and, from Propositions \ref{prop: laplacian and spherical potentials}, \ref{prop: properties of the laplacian}, we have

\[U_{\Delta \phi}(x) = \int \log \frac{1}{||x, y||_v} \; d(\Delta \phi)(y) = \int \phi(y)\;d\Delta_y \log\frac{1}{||x,y||} =  \phi(x) - \int \phi(y)\; d\lambda_{\Ar, v}(y)\]

so that $U_{\Delta \phi}$ differs from $\phi$ by a constant. In particular $U_{\Delta \phi}$ is continuous. By Proposition \ref{prop: laplacian and spherical potentials}, $\Delta U_{\mu_{n, v}} = \mu_{n, v} -  \mu_{n, v}(\Pberk)\lambda_{\Ar,v}$ and therefore 
\begin{align*} 
\int \phi \; d\mu_{n, v}  &= \int f\; d( \Delta U_{\mu_{n, v}} + \mu_{n, v} (\Pberk)\lambda_{\Ar, v})\\
                              &= \int U_{\mu_{n, v}}  \; d \Delta \phi + \mu_{n, v} (\Pberk) \int \phi\; d\lambda_{\Ar, v}\\
                              &= \< \mu_{n, v} , \Delta \phi \>_v + \mu_{n, v}(\Pberk) \int \phi\; d\lambda_{\Ar, v}.
\end{align*}
As $\log ||x, y||_v^{-1} \geq 0$, the local form $\< \cdot, \cdot\>_v$ is always at least positive semi-definite; in particular, the Cauchy-Schwarz inequality applies. Thus $|\< \mu_{n, v} , \Delta \phi \>_v| \leq \| \mu_{n, v} \|_v \| \Delta \phi\|_v$ and so 
\[  \bigg| \int \phi \; d \mu_{n, v} \bigg| \leq \| \mu_{n, v} \|_v \| \Delta \phi\|_v + \left| \mu_{n, v}(\Pberk)\right|\int |\phi|\; d\lambda_{\Ar, v}. \]
 From Corollary \ref{cor: pos def of ag inner product}, that $\|\mu_n\| \to 0$ implies $|\mu_{n,v}(\Pberk)| \to 0$ as $n \to \infty$, and by assumption $\|\mu_{n, v}\|_v^2$ also tends to $0$ as $n \to \infty$. Consequently $\int \phi\; d\mu_{n, v} \to 0 \qas n \to \infty$ for all $\phi \in E_v$, which shows the first statement. The second statement follows from the first and from the fact that $E_v$ is dense in $C(\Pberk)$.
\end{proof}

\begin{rem} 
It is easy to see that weak convergence at each place does not give convergence in the norm of $\Aad$. Here is a simple illustration. For each place $v \in K$, let $\nu_v$ be any log-continuous Radon measure of total mass $0$ with $r_v\|\nu_v\|_v^2 = 1$. Enumerating the places of $K$ by $v_n$ ($n \geq 1)$, let $\mu_n$ be the adelic measure which is the zero measure at all places $v \neq v_n$ and $\mu_{n, v_n} = \nu_v$. Then $\|\mu_n\|^2 = 1$ for all $n$, yet $\mu_{n, v} \to 0$ weakly at each place. 
\end{rem}

The following observation will be useful in the following section. 
\begin{prop}\label{prop: potential at infinity}
    Let $\mu = (\mu_v)_{v \in M_K} \in \Aad$ be an adelic probability measure with $(\!(\mu, \mu)\!) = 0$. Then $\|\mu\|^2 = 2\sum_{v \in M_K} r_v U_{\mu, v}(\infty)$. In particular, $\|\mu_f\|^2 = 2\sum_{v \in M_K} r_v U_{f, v}(\infty)$ whenever $f$ is a polynomial. 
\end{prop}

\begin{proof}
By (\ref{eq: chordal to affine dist}) we have
\begin{align*}
    \|\mu_v\|^2_v &= \iint \log \frac{1}{\delta_v(x,y)} \, d\mu_v(x)\,d\mu_v(y) + \iint \log \frac{1}{||x,\infty||_v} \, d\mu_v(x)\, d\mu_v(y) + \iint \log \frac{1}{||y, \infty||_v}\, d\mu_v(x)\, d\mu_v(y)\\
    &= (\!( \mu_v, \mu_v)\!)_v  + 2U_{\mu,v}(\infty).
\end{align*} 
Thus $\|\mu\|^2 = \sum_{v \in M_K} r_v ((\!( \mu_v, \mu_v)\!)_v  + 2U_{\mu,v}(\infty)) = (\!(\mu, \mu)\!) + 2\sum_{v \in M_K} r_v U_{\mu, v}(\infty) = 2\sum_{v \in M_K} r_v U_{\mu, v}(\infty)$. As remarked in the introduction, a result of Baker and Hsia implies $(\!(\mu_f, \mu_f)\!)  =0$ whenever $f$ is a polynomial \cite{baker2005canonical}; this also follows from Lemme 5.3 of \cite{favre2006equidistribution}. Therefore $\|\mu_f\|^2= 2\sum_{v \in M_K} r_v U_{f, v}(\infty)$ in this case. 
\end{proof}
\section{Comparison With the Arakelov Height: Monic Polynomials}\label{sec: comparison w ar height MONIC}

In both this section and the next, the pull-back of a measure on $\Pberk$ by a rational map $f$ will play an important role. We briefly summarize the discussion in chapter 9 of \cite{baker2010potential}. Given a rational map $f$ defined over $\C_v$ and a continuous function $\phi$ defined on $f^{-1}(V)$ of $\Pberk$, where $V$ is a subset of $\Pberk$, the \textit{push-forward} of $\phi$ by $f$ is the function  $(f_*\phi)(z) = \sum_{f(\zeta) = z} m_f(\zeta) \phi(\zeta)$. Here $m_f(\zeta)$ is an extension of the notion of algebraic multiplicity to $\Pberk$. In particular, if $z \in \C_v$ and $\alpha_i$ are the roots of $f(\zeta) = z$ counting multiplicity, then $(f_*\phi)(z) = \sum_{i=1}^{\deg(f)} \phi(\alpha_i)$. The function $f_*\phi$ is defined and continuous on $V$. Given  a Radon probability measure $\mu$ on $\Pberk$, the \textit{pull-back} of $\mu$ by $f$, written $f^* \mu$, is defined to be the unique Radon measure with  $\int \phi\; df^*\mu = \int f_*\phi\; d\mu$. From this definition and the fact that $\sum_{f(\zeta) = z}m_f(\zeta) = \deg(f)$ for all $z\in \Pberk$ we find that $f_* (f^*\mu) = \deg(f)\cdot \mu$, where $f_*(\nu)$ is the usual push-forward of a measure $\nu$. For proofs of these claims, see chapter 9 of \cite{baker2010potential}. We recall that the canonical measure of $f$ on $\Pberk$ is characterized among log-continuous probability measures by  $f^* \mu_{f, v} = \deg(f) \mu_{f, v}$.\\

The \textit{Arakelov height of a rational map} $f$ of degree $d$ is the Arakelov height of the coefficients of $f$, viewed as a point in projective space. Thus, if $f(z) = \frac{a_dz^d + \cdots + a_0}{b_d z^d + \cdots + b_0}$ (where one of $a_d$ or $b_d$ is non-zero) is defined over $K$, then  $h_{Ar}(f) = \sum_{v \in M_K} r_v \log |f|_v$, 
where
\[  
|f|_v = \begin{cases}
\max_{0\leq i, j \leq d} (|a_i|_v, |b_j|_v) & \qif v \nmid \infty,\\
\left( \sum_{i=0}^d |a_i|_v^2 + |b_i|_v^2\right)^{1/2} & \qif v \mid \infty.
\end{cases}
\]

In particular, if $f(z) = z^d + a_{d-1}z^{d-1} + \cdots + a_0$ is a monic polynomial defined over $K$, then $|f|_v = \left(2 + \sum_{i=0}^{d-1} |a_i|_v^2\right)^{1/2}$ for Archimedean $v$ and $|f|_v = \max_i (1, |a_i|_v)$ at finite places $v$ of $K$. The main theorem of this section is the following estimate on the norm $\|\mu_f\|$ when $f$ is a monic polynomial.

\begin{thm}\label{thm: norm of monic poly is a height}
Let $f(z) = z^d + a_{d-1}z^{d-1} + \cdots + a_0$ be a monic polynomial of degree $d \geq 2$ defined over $K$. Then $\frac{d}{2}\|\mu_f\|^2 = h_{\Ar}(f) + O(1)$, where the implied constants depend only on $d$.
\end{thm}

This and the second item of Proposition \ref{prop: ag inner product essentially standard inner product} imply Theorem \ref{thm: heigh on monic poly} from the introduction.  The strategy for proving Theorem \ref{thm: norm of monic poly is a height} is to make use of the \textit{(classical) filled Julia set}, which, following \cite{benedetto2019dynamics}, we denote by $\Kcal_{I, v}$. This is defined to be those elements of $\C_v$ which are bounded under iteration by $f$, hence $\Kcal_{I, v} = \{ z \in \C_v : \lim_{n\to \infty} |f^n(z)|_v \neq \infty\}$. When $v$ is Archimedean, this is the usual filled Julia set from complex dynamics. For finite places $v$, this consists of the Type I points of the \textit{Berkovich filled Julia set}, $\Kcal_v$,  which is defined by

\[  \Kcal_v = \{ \zeta \in \Pberk : \lim_{n \to \infty} f^n(\zeta) \neq \infty\},\; \text{so that}\; \Kcal_{I, v} = \Kcal_v \cap \C_v. \]

At each place $v$, the \textit{Julia set} $\Jcal_v$ can be defined as the topological boundary of $\Kcal_v$. As in the classical case, $\mu_{f, v}$ is supported on $\Jcal_v$. We can therefore control the norm $\|\mu_f\|^2$ via escape rate arguments at each place. Namely, we employ escape rate arguments to prove the following. 

\begin{lem}\label{lem: f_y vs f escape rate arguments}
Let $f$ be a monic polynomial of degree $d\geq 2$ defined over $\C_v$ and let $y \in \Kcal_{I, v}$. Let $f_y(z)$ be the polynomial $f(z) - y$. There are positive constants $A_v$ and $B_v$ depending only on the degree $d$ so that $A_v |f|_v \leq  |f_y|_v \leq B_v |f|_v$. Furthermore, $A_v = B_v = 1$ when $v$ is non-Archimedean.
\end{lem}

We defer the proof of Lemma \ref{lem: f_y vs f escape rate arguments} to the end of this section. We will also need to know how the roots of a monic polynomial are related to the quantities $|f|_v$. 

\begin{lem}\label{lem: monic polynomial roots and coefficients}
Let $f$ be a monic polynomial of degree $d\geq 2$ defined over $\C_v$. Factor $f(z)$ over $\C_v$ as $f(z) = \prod_{i=1}^d (z-\alpha_i)$. There are positive constants $C_v$ and $D_v$ depending only on the degree $d$ so that  $C_v |f|_v \leq \prod_{i=1}^d ||\alpha_i, \infty||_v^{-1} \leq D_v |f|_v$. Furthermore, $C_v = D_v = 1$ when $v$ is non-Archimedean 
\end{lem}

This is proved as part of Theorem VIII.5.9 of \cite{silverman2009arithmetic}. (Silverman uses $\max(1, |\alpha|_v)$ at the infinite places, but, with minor modifications, the same proof works with $||\alpha_i, \infty||_v^{-1} = \sqrt{1 + |\alpha|_v^2}$ instead of $\max(1, |\alpha|_v)$.) We turn to the proof of Theorem \ref{thm: norm of monic poly is a height}.

\begin{proof}[Proof of Theorem \ref{thm: norm of monic poly is a height}]
By Proposition \ref{prop: potential at infinity}, we have $\sum_{v \in M_K} r_v U_{\mu, v}(\infty) = \frac{1}{2}\|\mu_f\|^2$. The characterizing equation $f^*\mu_{f, v} = d\mu_{f, v}$ implies 
\begin{equation}\label{eq: norm of monic poly and pot} 
\frac{d}{2}\|\mu_f\|^2  = d\sum_{v\in M_K} r_v U_{f, v}(\infty) =   \sum_{v \in M_K} r_v \int \sum_{f(z) = y } m_f(z)\log \frac{1}{||z, \infty||_v} \;d\mu_{f, v}(y). 
\end{equation}

Thus, it suffices to bound $dU_{f, v}(\infty) =  \int \sum_{f(z) = y } m_f(z)\log ||z, \infty||^{-1}_v\;d\mu_{f, v}(y)$ in terms of $|f|_v$.\\

Suppose that $v$ is non-Archimedean. Define $\psi: \Aberk \to [0, \infty)$ by $\psi(\zeta) = -\log ||\zeta, \infty||_v$. Then $\psi$ is continuous on $\Aberk$. That $f$ fixes the point at infinity implies that the push forward $(f_*\psi)(\zeta) = \sum_{f(\alpha) = \zeta} m_f(\alpha) \log||\alpha, \infty||_v^{-1}$ is also continuous on $\Aberk$. Let $y \in \Kcal_{I, v}$. Factor $f(z) - y$ over $\C_v$ as $f(z) - y = \prod_{i=1}^d (z- \alpha_i)$, counting multiplicities. Then $(f_* \psi)(y) = \log \prod_{i=1}^d ||\alpha_i, \infty||_v^{-1}$.  Apply Lemma \ref{lem: monic polynomial roots and coefficients} to the polynomial $f_y(z) = f(z) - y$ to conclude that $(f_* \psi)(y) = \log |f_y|_v$. From Lemma \ref{lem: f_y vs f escape rate arguments}, then $(f_*\psi)(y) = \log |f_y|_v = \log |f|_v.$ That $\C_v$ is dense in $\Aberk$ implies that $\Kcal_{I, v}$ is dense in $\Kcal_v$. Thus $f_*\psi$ is a continuous function which is equal to the constant $\log |f|_v$ on a dense subset of $\Kcal_v$. Hence $(f_*\psi)(\zeta) = \log |f|_v$ for all $\zeta \in \Kcal_v$. As $\supp\; \mu_{f, v} \subset \Kcal_v$, then  $dU_{f, v}(\infty) = \int (f_*\psi)\; d\mu_{f, v} = \log |f|_v$.\\

The Archimedean case is similar: let $y \in \Kcal_v$ and factor $f(z) - y$ over $\C$ as $f(z) - y = \prod_{i=1}^d (z - \alpha_i)$, counting multiplicities. By Lemmas \ref{lem: monic polynomial roots and coefficients} and \ref{lem: f_y vs f escape rate arguments} we find that $\log \prod_{i=1}^d \sqrt{1 + |\alpha_i|_v^2} = \log |f_y|_v + O(1)$ and that $\log |f_y|_v = \log |f|_v + O(1)$; therefore $\log \prod_{i=1}^d \sqrt{1 + |\alpha_i|_v^2} = \log |f|_v + O(1)$.
As $\supp\; \mu_{f, v} \subset \Kcal_v$, we simply integrate this with respect to $y$ to find $dU_{f, v}(\infty) = \log |f|_v + O(1)$,
where, as in the statements of Lemmas \ref{lem: monic polynomial roots and coefficients} and \ref{lem: f_y vs f escape rate arguments}, the implied constants depend only on the degree of $f$. Multiplying $dU_{f, v}(\infty)$ by $r_v$, summing over all places, and applying (\ref{eq: norm of monic poly and pot}) gives the theorem.
\end{proof}

\begin{proof}[Proof of Lemma \ref{lem: f_y vs f escape rate arguments}]

Suppose $v \nmid \infty$. An elementary escape rate argument shows that $|y|_v \leq |f|_v$ for all $y\in \Kcal_{I, v}$. Thus, if $|a_k|_v = |f|_v =  \max_{0 \leq j\leq d-1}(1, |a_j|_v)$ for some $k > 0$, then $|a_k|_v \geq |a_0|_v$ and $|a_k|_v \geq |y|_v$. So $|a_k|_v \geq \max (|a_0|_v, |y|_v ) \geq |a_0 - y|_v$. Hence 
\[|f_y|_v = \max(1, |a_j|_v, |a_0 - y|_v) = |f|_v.\]
Similarly, if $|f|_v= 1$ then $|f_y|_v = 1$. So assume $|a_0|_v = |f|_v > 1$ and  $|a_k|_v < |f|_v$ for $1\leq k \leq d-1$. We claim that then $|y|_v < |f|_v$. Indeed, if $|y|_v = |f|_v$ then $|y^d|_v > |a_j y^j|_v$ for $0\leq j < d$ so that $|f(y)|_v = |y|^d$ (by the equality case of the strong triangle inequality). Iterating this equation shows that $|f^n(y)|_v \to \infty$ so that $y \notin \Kcal_{I, v}$. Thus $|y|_v < |f|_v$ and so $|a_0 - y|_v = |a_0|_v = |f|_v$; it follows that $\max_{1\leq j\leq d-1}( 1, |a_j|_v, |a_0 - y|_v) = \max_{0\leq j\leq d-1}(1, |a_j|_v)$. \\

We turn to the Archimedean case. To ease notation, we omit the dependence on $v$. The quantity $E_f := 1 + \sum_{j=0}^{d-1} |a_j|$ is more amenable to study than $|f|$. Clearly $(1/2)|f| \leq E_f$, and, by Jensen's inequality applied to the map $t \mapsto t^2$, we find that $E_f \leq \sqrt{d}|f|_v$. So it is enough to show that the lemma holds at Archimedean places with $E_f$ in place of $|f|$. \\

An easy escape rate argument shows that $|y| \leq E_f$ for $y$ in the filled Julia set. By the triangle inequality
\[E_{f_y} = 1 + |a_0 - y| +  \sum_{j=1}^{d-1} |a_j| \leq 1 + |a_0| + E_f + \sum_{j=1}^{d-1} |a_j| \leq 2E_f,\]
which gives the upper bound.  For the lower bound, we break into cases depending on the size of the constant term. Suppose that $1 + \sum_{j=1}^{d-1} |a_j| \geq \epsilon |a_0|$, where $\epsilon \in (0,1/2)$ is an absolute constant to be determined later. Thus $1 + \sum_{j=1}^{d-1} |a_j| \geq \frac{\epsilon}{1+\epsilon}E_f$, which implies
\[  E_{f_y} = 1 + |a_0 - y| +  \sum_{j=1}^{d-1} |a_j| \geq \frac{\epsilon}{1 + \epsilon}E_f.  \]

So assume that $1 + \sum_{j=1}^{d-1} |a_j| < \epsilon |a_0|$. We claim that then $|y| \leq (1-\epsilon)|a_0|$. Suppose to the contrary that $|y| > (1-\epsilon)|a_0|$. Then $|y| > \epsilon^{-1}(1-\epsilon) > 1$, as  $|a_0| > \epsilon^{-1}$ and $\epsilon < 1/2$.  Now
\begin{align*}
|f(y)| &= |y^d + a_{d-1}y^{d-1} + \cdots + a_0| \\
       &\geq |y|^d  - \left(|a_0| + |y|^{d-1} \sum_{j=1}^{d-1}|a_j|\right)\\
       &\geq |y|^d - \left(|a_0| + |y|^{d-1}(\epsilon|a_0| - 1)\right)\\
       &=|y|^d + |y|^{d-1} - |a_0| - \epsilon|a_0| |y|^{d-1}.
\end{align*}

That $|y|\epsilon(1-\epsilon)^{-1} > 1$ implies
\[ \frac{|y|^d}{1 + \epsilon|y|^{d-1} }\geq \frac{|y|^d}{\frac{\epsilon}{1-\epsilon}|y|^{d-1} + \epsilon|y|^{d-1}  } = |y|\frac{1-\epsilon}{\epsilon(2-\epsilon)}. \]
Now choose $\epsilon$ to be the smallest real root of $(1-T)^2 - T(2- T) =  0$; i.e., $\epsilon = 1 - 1/\sqrt{2}$. Then $\epsilon \in (0,1/2)$ and $(1-\epsilon)(\epsilon(2 - \epsilon))^{-1} =(1-\epsilon)^{-1}$ so that 
\[ \frac{|y|^d}{1 + \epsilon|y|^{d-1} } \geq |y|\frac{1-\epsilon}{\epsilon(2-\epsilon)} = \frac{|y|}{1-\epsilon} > |a_0|, \]
where the final inequality is by assumption on $|y|$. Therefore
\[ |f(y)| \geq  |y|^d + |y|^{d-1} - |a_0| - \epsilon|a_0| |y|^{d-1} > |y|^{d-1}.  \]
Iterating this inequality shows that $|f^n(y)| \to \infty$ as $n \to \infty$. But this contradicts $y$ being in the filled Julia set. So $|y| \leq (1-\epsilon)|a_0|$ in this case, hence $|a_0 - y| \geq |a_0| - (1-\epsilon)|a_0| = \epsilon|a_0|$. Consequently
\[ E_{f_y} \geq 1 + \epsilon|a_0| + \sum_{j=1}^{d-1} |a_j| \geq \epsilon E_f \geq \frac{\epsilon}{1 + \epsilon}E_f.  \]
Thus, regardless of the size of the constant term, we have $E_{f_y} \geq \frac{\epsilon}{1 + \epsilon} E_f$. 

\end{proof}

\section{Comparison with the Arakelov Height: Rational Maps}\label{sec: comparison w ar ht rat}
We turn to our main theorem (Theorem \ref{thm: height on polynomials}). 
\subsection{An Adjoint Formula}
We define the pull-back and push-forward of an adelic measure \\$\mu = (\mu_v)_{v \in M_K} \in \Aad$ coordinate-wise: $f^*\mu = (f^*\mu_v)_{v \in M_K}$ and $f_*\mu = (f_* \mu_v)_{v \in M_K}$. It is useful to know how the pull-back and push-forward interact with our inner product.  For measures of total mass zero, the push-forward and pull-back are  adjoints of each other. 

\begin{prop}\label{prop: Adjoint}
Suppose that $\mu$ and $\nu$ are in $\Aad^0$. If $f$ is a non-constant rational map defined over $K$ then  $\< f^* \mu, \nu\> = \< \mu, f_* \nu\>$.
In particular, $\| f^* \mu\|^2 = d\|\mu\|^2$.
\end{prop}

The proof is entirely local, and indeed the adjoint relationship holds for local forms as well. The key is the following relationship between the Laplacian and the action of a rational map: if $\varphi\in  \BDV_v(\Pberk)$ then 
\begin{equation}\label{eq: laplacian adjoint}  f^* \Delta(\varphi) = \Delta( \varphi \circ f). \end{equation}
See Proposition 9.56 of \cite{baker2010potential} for a proof, or the remarks in Section 6.1 of \cite{favre2006equidistribution}. 

\begin{proof}[Proof of Proposition \ref{prop: Adjoint}]
The proof is essentially a reformulation of (\ref{eq: laplacian adjoint}). Since $U_{\mu, v}$ is in $\BDV_v$ for all $v$,  and since $\Delta U_{\mu, v} = \mu_v$, $\Delta U_{\nu, v} = \nu_v$ as $\mu_v$ and $\nu_v$ have total mass $0$, then 
\[    \< f^* \mu_v, \nu_v\>_v = \int U_{\nu, v}\; df^* \mu_v
                               = \int U_{\nu, v}d \Delta( U_{\mu, v} \circ f)
                               = \int U_{\mu, v} \circ f\; d\nu_v 
                               = \< \mu_v, f_* \nu_v\>_v. \]
The second statement follows from the first: $\< f^* \mu_v, f^* \mu_v\>_v = \< \mu_v, f_*f^* \mu_v\>_v = d\<\mu_v, \mu_v\>_v$. Multiplying by $r_v$ and summing over all places gives the proposition.
\end{proof}

\subsection{Comparison with the pull-back of the Arakelov Measure.}

The Arakelov height of a given rational map $f$ in some sense measures the arithmetic complexity of $f$. It is therefore natural to wonder if the norm $\|\mu_f\|^2$ of $f$, which is a measure of arithmetic-dynamical size, can be related to Arakelov height of $f$. It turns out that $\|\mu_f\|^2$ can be bounded in terms of the Arakelov height of $f$. 

\begin{thm}\label{thm: arakelov height norm}
Let $f$ be a rational map of degree $d \geq 2$, defined over a number field $K$. There are positive constants $c_1, c_2, c_3, c_4$ depending only on the degree of $f$  so that
\[ c_1 h_{Ar}(f) - c_2 \leq \frac{d}{2}\|\mu_f\|^2 \leq c_3 h_{Ar}(f) + c_4.  \]
Moreover $c_1, c_3 \to 1$ as $d \to \infty$ and $c_2, c_4$ grow linearly in $d$. 
\end{thm}

As $h_\Ar = h + O(1)$ for any Weil height $h$ on $\Rat_d(K)$, with the implied constants depending on $d$ and $h$, Theorem \ref{thm: arakelov height norm} and Proposition \ref{prop: ag inner product essentially standard inner product} imply Theorem \ref{thm: height on polynomials} from the introduction. When examining the norm associated to monic polynomials, we could deduce information on $\mu_{f, v}$ via the filled Julia sets $\Kcal_{I, v}$ at each place. In the case of rational maps, such tools are not available and we are forced to take a more circuitous route in proving Theorem \ref{thm: arakelov height norm}. Rather than working with the norm of $\mu_f$ directly, we instead study the pull-back of the Arakelov adelic measure $\lambda_{Ar}$ by $f$. As the following simple lemma shows, the norms of these two adelic measures are never too far apart.

\begin{lem}\label{lem: pull back bounds}
Let $\mu_f = (\mu_{f, v})_{v \in M_K}$ be the canonical adelic measure of a rational map $f$ defined over $K$ with degree $d \geq 2$. Then  

\begin{equation}\label{eq: pull back norm inequality}  
\frac{d}{(\sqrt{d} + 1)^2} \|d^{-1}f^* \lambda_{\Ar}\|^2 \leq \|\mu_f\|^2 \leq\frac{d}{(\sqrt{d} -1)^2} \|d^{-1}f^*\lambda_\Ar\|^2.
\end{equation}
\end{lem}

\begin{proof}
    From Proposition \ref{prop: Adjoint} and the fact that $d^{-1}f^* \mu_f = \mu_f$, we have
    \[  \| d^{-1} f^*\lambda_\Ar - \mu_f\|^2 = \| d^{-1}f^*( \lambda_\Ar - \mu_f)\|^2 = d^{-1} \|\lambda_\Ar - \mu_f\|^2.   \]
    Now $\|\lambda_\Ar - \mu_f\|^2 = \|\mu_f\|^2 - 1/2$, by Proposition \ref{prop: ag inner product essentially standard inner product}. (Recall that $\< \cdot, \cdot\> = \< \cdot, \cdot\>_{\lambda_\Ar, 1/2}$.) So the triangle inequality gives
    \[ \bigg|\| d^{-1} f^*\lambda_{\Ar}\| - \|\mu_f\| \bigg| \leq \|d^{-1}f^*\lambda_\Ar - \mu_f\| \leq  d^{-1/2} \sqrt{ \|\mu_f\|^2 - \frac{1}{2}} \leq d^{-1/2} \|\mu_f\|. \]
    Simple manipulation of this inequality shows that   
    \[  \frac{\sqrt{d}}{\sqrt{d} + 1}\|d^{-1}f^*\lambda_\Ar\| \leq\|\mu_f\| \leq \frac{\sqrt{d}}{\sqrt{d} -1}\|d^{-1} f^*\lambda_{\Ar}\|. \]
    Squaring these inequalities gives (\ref{eq: pull back norm inequality}).

    
       
\end{proof}

So, to prove Theorem \ref{thm: arakelov height norm}, it suffices to bound  $\|f^* \lambda_{\Ar}\|^2$ in terms of the Arakelov height of $f$. 

\subsection{Bounds on the pull-back of the Arakelov Measure}\label{sec: bounds on pull back}

We are going to prove a fact which may be of interest in its own right: for $f \in \Rat_d(K)$,  $(d/2)\|d^{-1}f^*\lambda_{\Ar}\|^2 = h_\Ar(f) + O(1)$, with implied constants depending only on $d$. The strategy is to compare $\|f^* \lambda_{\Ar, v}\|^2_v$ with $|f|_v$ at each place of $K$. At  finite places, we obtain an exact formula (Lemma \ref{lem: non-arch pull back of sigma}), and at Archimedean places we obtain an integral representation which is easily bounded in terms of the Archimedean contribution to $h_\Ar (f)$ (Propositions \ref{prop: arch pull back of sigma} and \ref{prop: bounds on archimedean pull back}). We begin by fixing some notation.  Write 

\[f(z) = \frac{P(z)}{Q(z)} = \frac{a_d z^d + \cdots + a_0}{b_d z^d + \cdots + b_0},\]
where at least one of $a_d$ or $b_d$ is non-zero and where $P$ and $Q$ have no common zero. Set $p = \deg P$ and $q = \deg Q$, and $d = \deg f =  \max(p,q)$. So $a_p$ is the leading coefficient of $P(x)$ and $b_q$ is the leading coefficient of $Q(x)$. For non-Archimedean places $v$ we write $|\cdot|_v$ for the continuous extension of the $v$-adic absolute value on $\C_v$ to $\Aberk$. Explicitly, $|x|_v = \delta_v(x, 0)$, where $\delta_v$ is the Hsia kernel from Section 2. Equivalently, $|x|_v = [T]_x$, where $[\cdot]_x \in \Aberk$ is viewed as a seminorm on $\C_v[T]$ (see chapters 1 and 2 of \cite{baker2010potential} for the semi-norm definition of Berkovich spaces). The following consequence of Jensen's formula and its non-Archimedean analogue will be central to what follows. 

\begin{prop}\label{prop: arakelov jensen}
Let $v \in M_K$ and let $x,y \in \C_v$. Then
\[ \int \log | x - z y|_v \; d\lambda_{\Ar, v}(z) = 
\begin{cases}
\log\max (|x|_v, |y|_v) &\qif v \nmid \infty;\\
\log \sqrt{ |x|^2_v + |y|^2_v } & \qif v \mid \infty.
\end{cases}
\]
\end{prop}

\begin{proof}
The proposition is trivial if $y =0$, so assume  $y$ is non-zero. Suppose $v$ is non-Archimedean.  Then $\log [T- (x/y)]_{\zeta_\Gauss} = \log \sup_{|z|_v\leq 1} |z - (x/y)|_v = \log^+ |x/y|_v$, and so  
\[\int \log |x - z y|_v \; d\lambda_{\Ar, v}(z) = \log|y|_v + \int \log |(x/y) - z|_v\; d\lambda_{\Ar, v}(z) = \log \max (|x|_v, |y|_v).\]
Suppose that $v$ is Archimedean and that $y = 1$. We write $\lambda_{\Ar, v}$ in polar coordinates, suppress the dependence on $v$, and apply Jensen's formula:
\begin{align*}
        \int \log|x - z|\; d\lambda_{\Ar}(z) &= \int_0^\infty    \left( \int_0^{2\pi} \log|x - re^{it}|\frac{dt}{2\pi} \right) \frac{2r \ dr}{(1 + r^2)^2}\\
                                  &= \int_0^\infty \log \max( r, |x|) \; \frac{2\,dr}{(1 + r^2)^2} \;\;\;\;\;\; (\text{by Jensen's formula})\\
                                  &= \log |x| \int_{0}^{|x|} \frac{2r dr}{(1 + r^2)^2}  + \int_{|x|}^\infty \log r \frac{2 r dr}{(1 + r^2)^2}\\
                                  &= \frac{ |x|^2 }{1 + |x|^2}\log (|x|)   + \frac{1}{2}\log( 1 + |x|^{-2}) + \frac{\log |x|}{1 + |x|^2}\\
                                  &= \log\sqrt{1 + |x|^2}.
\end{align*}
Applying the above computation to $\log |x - yz| = \log |y| + \log |(x/y) -z |$ finishes the proof.  
\end{proof}

We define a function $I_v: \C_v \to \R$ by 

\[  I_v(w) = \int \log |P(z)-wQ(z)|_v\; d\lambda_{\Ar, v}(z). \]
For non-Archimedean places $v$, $I_v(w) = \log |(P- wQ)(\zeta_\Gauss)|_v = \log \sup_{|z|_v \leq 1} |P(z)-wQ(z)|_v$. Of course, $I_v$ depends on how $P$ and $Q$ are normalized; as we will summing over all places in the end, this will be irrelevant due to  product formula. That $\lambda_{\Ar, v}$ is log-continuous implies $I_v$ exists for any $w\in \C_v$. With the preliminaries dealt with, we turn to the following key lemma, which relates the pull-back by $f$ with our potential kernel $\log ||x, y||_v^{-1}$. We phrase this lemma so that it holds at all places.

\begin{lem}\label{lem: resultant pull back}

Write $f(z) = \frac{P(z)}{Q(z)}$ as above. Fix elements $x$ and $y$ of  $\C_v$ with $x \neq y$, and with $x, y$ not equal to $f(\infty)$.  For $1 \leq i, j \leq d$, let $\alpha_i$ and $\beta_j$ be respectively the roots of $f(z) = x$ and $f(z) = y$ as polynomials in $z$, counting multiplicities. Then 
\[ \sum_{i,j =1}^d \log \frac{1}{||\alpha_i, \beta_j||_v} =  -\log|a_p^{d-q} b_q^{d-p}\Res(P, Q)|_v  - d\log|x-y|_v +  dI_v(x) + dI_v(y), \]
\end{lem}


\begin{proof}
The proof is a straightforward, albeit rather tedious, computation. Note $|\Res(P, Q)|_v \neq 0$ as $P$ and $Q$ have no common zero. That $x \neq y$ implies both $||\alpha_i , \beta_j||_v$ and $|x - y|_v$ are not zero. Moreover, the polynomials $P(z) - xQ(z)$ and $P(z) - yQ(z)$ both have degree $d$; for if not, then either $x$ or $y$ would equal $f(\infty)$.  Factor $P(z) - x Q(z)$ and $P(z) - y Q(z)$ over $\C_v$ as 
\[P(z) - x Q(z) = \ell_x \prod_{i=1}^d (z - \alpha_i) \qand P(z) - y Q(z) = \ell_y\prod_{j=1}^d (z - \beta_j)\]
where $\ell_x = a_d - x b_d$ and $\ell_y = a_d - y b_d$. By definition of $\beta_j$,  

\[\prod_{i,j=1}^d |\alpha_i - \beta_j|_v =  \prod_{i=1}^d |\ell_y|^{-1}_v|P(\alpha_i) - y Q(\alpha_i)|_v, \]

and, as $P(\alpha_i) = x Q(\alpha_i)$,  

\[\prod_{i,j=1}^d |\alpha_i - \beta_j|_v = |\ell_y|^{-d}_v |x-y|_v^d \prod_{i=1}^d |Q(\alpha_i)|_v.\]

Factor $Q(z)$ over $\C_v$ as $Q(z) = b_q \prod_{k=1}^q (z - r_k)$. Then 

\begin{align*}
    \prod_{i=1}^d |Q(\alpha_i)|_v &= \prod_{i=1}^d |b_q|_v \prod_{k=1}^q |\alpha_i - r_k|_v  \\
                             &= |b_q|_v^{d} \prod_{k=1}^q \prod_{i=1}^d |\alpha_i - r_k|_v\\
                             &= |b_q|_v^{d} \prod_{k=1}^q |\ell_x|_v^{-1} |P(r_k) - x Q(r_k)|_v\\
                             &= |b_q|_v^{d} |\ell_x|_v^{-q} \prod_{k=1}^q |P(r_k)|_v\\
                             &= |\ell_x|_v^{-q} |b_q|_v^{d- p} |\Res(P, Q)|_v.
\end{align*}  

Thus 

\begin{equation}\label{eq: 1 lem pull back resultant}
    \prod_{i,j=1}^d |\alpha_i - \beta_j|_v = |\ell_y|_v^{-d} |\ell_x|_v^{-q} |b_q|_v^{d-p} |x-y|_v^d |\Res(P, Q)|_v.
\end{equation}  

Now assume that $v \nmid \infty$. By Proposition \ref{prop: arakelov jensen}

\[ \log |(P-xQ)(\zeta_\Gauss)|_v = \int\log \left\{|\ell_x|_v \prod_{i=1}^d |z - \alpha_i|_v\right\} d\lambda_{\Ar, v} = \log |\ell_x|_v + \sum_{i=1}^d \log^+ |\alpha_i|_v  \]




Similarly, $\log|(P-yQ)(\zeta_\Gauss)|_v = \log \left( |\ell_y|\prod_{j=1}^d \max(1, |\beta_j|_v) \right)$. As $\alpha_i$ and $\beta_j$ are in $\C_v$, then combining (\ref{eq: 1 lem pull back resultant}) with the expression for $||\alpha_i, \beta_j||_v$ in affine coordinates yields 

\begin{align*} 
    \sum_{i,j=1}^d \log \frac{1}{||\alpha_i, \beta_j||_v} &= -\sum_{i,j=1}^d \log  \frac{|\alpha_i - \beta_j|_v}{\max(1, |\alpha_i|_v)\max(1, |\beta_j|_v)} \\
                                     &= -\log \prod_{i, j = 1}^d |\alpha_i - \beta_j|_v  + d\log \prod_{i=1}^d \max(1, |\alpha_i|_v) + d\log \prod_{j=1}^d \max(1, |\beta_j|_v)\\
                                     &= -\log\left( |\ell_y|_v^{-d} |\ell_x|_v^{-q} |b_q|_v^{d-p} |x-y|^d_v |\Res(P, Q)|_v\right) +\\
                                     &\;\;\; + (dI_v(x) - d\log |\ell_x|_v) + (d I_v(y) - d\log |\ell_y|_v)\\
                                     &= -\log C_v(x, y)  - d\log |x-y| + dI_v(x) + dI_v(y)
\end{align*}
where $C_v(x,y) = |\ell_x|_v^{d-q} |b_q|_v^{d-p} |\Res(P,Q)|_v$. Now if $d = q$ then $|\ell_x|^{d-q} = 1 = |a_p|^{d-q}$. If $d > q$ then $d = p$ and $b_d = 0$ so that $\ell_x = a_d - b_d x = a_d$. Thus $C_v(x, y) = |a_p|_v^{d-q} |b_q|_v^{d-q}|\Res(P, Q)|_v$, which proves the lemma for non-Archimedean $v$.  The argument for Archimedean places $v$ is similar.  From Proposition \ref{prop: arakelov jensen}

\begin{align*}
    \log \prod_{i=1}^d \sqrt{1 + |\alpha_i|_v^2} &= \int \log \prod_{i=1}^d |\alpha_i - w|_v \; d\lambda_{\Ar, v}(w)\\
                                            &= \int\log |\ell_x|_v^{-1} |P(w) - x Q(w)|_v \; d\lambda_{\Ar, v}(w)\\
                                            &= -\log |\ell_x|_v + I_v(x);
\end{align*} 
likewise $\log \prod_{j=1}^d \sqrt{1 + |\beta_j|_v^2} = -\log |\ell_y|_v + I_v(y)$. Combining this with Equation \ref{eq: 1 lem pull back resultant} shows that
\begin{align*}
    \sum_{i,j=1}^d \log \frac{1}{||\alpha_i, \beta_j||_v} &= -\log \prod_{i,j = 1}^d \frac{|\alpha_i - \beta_j|_v}{\sqrt{1 + |\alpha_i|_v^2}\sqrt{1 + |\beta_j|_v^2}}\\
                                      &= -\log \prod_{i,j = 1}^d |\alpha_i - \beta_j|_v + d\log \prod_{i=1}^d \sqrt{1 + |\alpha_i|_v^2} + d\log \prod_{j=1}^d \sqrt{1 + |\beta_j|_v^2} \\
                                     &= -\log\left( |\ell_y|_v^{-d} |\ell_x|_v^{-q} |b_q|_v^{d-p} |x-y|^d_v |\Res(P, Q)|_v\right) +\\
                                     &\;\;\; + (dI_v(x) - d\log |\ell_x|_v) + (d I_v(y) - d\log |\ell_y|_v)\\
                                     &= -\log (|a_p|_v^{d-q}|b_q|_v^{d-p}|\Res(P, Q)|_v) - d\log |x-y|_v + dI_v(x) + dI_v(y),
\end{align*}
where, as in the non-Archimedean case, the terms involving $\ell_y$ cancel and $|\ell_x|_v^{d-q} = |a_p|_v^{d-p}$.
\end{proof}


\subsection{Non-Archimedean Considerations}

At finite places $v$, we have an explicit relationship between $\log |f|_v = \log \max_{0\leq i,j\leq d} (|a_i|_v, |b_j|_v)$ and $\|f^*\lambda_{\Ar, v}\|^2_v$. 

\begin{lem}\label{lem: non-arch pull back of sigma}
For all $v\in M_K^0$, $\|f^*\delta_{\zeta_\Gauss}\|_v^2 = -\log |a_p^{d-q} b_q^{d-p}\Res(P,Q)|_v + 2d\log |f|_v$.
\end{lem}

Conceptually, the strategy for proving Lemma \ref{lem: non-arch pull back of sigma} is straightforward: we are going to take the limit as $(x,y) \to (\zeta_\Gauss, \zeta_\Gauss)$ of the expression in Lemma \ref{lem: resultant pull back}. However, because $\log ||x, y||_v^{-1}$ is not continuous on $\Pberk \times \Pberk$, some care is needed in how we carry out this limiting process. We start with a simple proposition. We recall that the Berkovich projective line over the $v$-adic complex numbers is metrizable (Corollary 1.20 of \cite{baker2010potential}). Thus we may work with sequences rather than nets. 

\begin{prop}\label{prop: non arch integral limit}
Suppose that $x_n \in \C_v$ with $x_n \to \zeta_\Gauss$. Then $I_v(x_n) \to \log |f|_v$.
\end{prop}

\begin{proof}
It follows from Proposition \ref{prop: arakelov jensen} that $\lim_{x_n \to \zeta_\Gauss} |a_i - x_n b_i|_v = \max(|a_i|_v, |b_i|_v)$. Therefore

\begin{align*}
    \lim_{x_n \to \zeta_\Gauss} I_v(x_n) &= \lim_{x_n \to \zeta_\Gauss} \log |(P-x_nQ)(\zeta_\Gauss)|_v \\
                                         &= \lim_{x_n \to \zeta_\Gauss} \log \max_{0 \leq i \leq d} |a_i - x_n b_i|_v\\
                                         &= \log \max_{0\leq i\leq d} \lim_{x_n \to \zeta_\Gauss} |a_i - x_n b_i|_v\\
                                         &= \log \max_{0\leq i\leq d} \max( |a_i|_v, |b_i|_v)\\
                                         &= \log |f|_v.
\end{align*} 
\end{proof}

\begin{proof}[Proof of Lemma \ref{lem: non-arch pull back of sigma}]
Fix $y \in \C_v$ and let $\beta_1, ..., \beta_d$ be the roots of $f(z) = y$, including multiplicities. Define a function $\psi_y: \Pberk \to \R \cup \{\infty\}$ by $\psi_y(z) = \sum_{j=1}^d \log||z, \beta_j||_v^{-1}$. We remark that $\psi_y$ is continuous everywhere except at the finitely many points $\beta_j$. Let $\{x_n\}_{n=1}^\infty$ be a sequence in $\C_v$ with $x_n \to \zeta_\Gauss$. By replacing $\{x_n\}_{n=1}^\infty$ with a subsequence if necessary, we assume that the conditions of Lemma \ref{lem: resultant pull back} are satisfied. In particular, $x_n \neq y$ and $\zeta_\Gauss \neq y$ for all $n$. Let $V$ be a sufficiently small open set containing $y$ so that $\zeta_\Gauss \notin V$ and $x_n \notin V$ for all $n$. As $\beta_j \in f^{-1}(V)$ for $1 \leq j \leq d$, then $\psi_y$ is continuous on $f^{-1}(\Pberk \setminus V)$. Therefore the push-forward of $\psi_y$ by $f$, $(f_* \psi_y)$, is continuous on $\Pberk \setminus V$. Since $x_n \to \zeta_\Gauss$ and  $x_n, \zeta_\Gauss \in \Pberk\setminus V$, then $(f_* \psi_y)(x_n) \to (f_*\psi_y)(\zeta_\Gauss)$. For $1 \leq i \leq d$ let $\alpha_{i, n}$  be the roots of $f(z) = x_n$, counting multiplicities. By Lemma \ref{lem: resultant pull back} then
\[  (f_* \psi_y)(x_n) = \sum_{1\leq i, j\leq d} \log\frac{1}{||\alpha_{i,n}, \beta_j||_v} = -\log |a_p^{d-q}b_q^{d-p}\Res(P,Q)|_v - d\log|x_n - y|_v + d I_v(x_n) + d I_v(y), \]
so that 
\begin{align*}
    (f_* \psi_y)(\zeta_\Gauss) &= \lim_{x_n \to \zeta_\Gauss} \left\{-\log |a_p^{d-q}b_q^{d-p}\Res(P,Q)|_v - d\log|x_n - y|_v + dI_v(x_n) + dI_v(y) \right\}\\
                            &= -\log|a_p^{d-q}b_q^{d-p}\Res(p,q)|_v + d\log|f|_v + d\log^+ |y|_v + dI_v(y),
\end{align*}  
where we have applied Proposition \ref{prop: non arch integral limit}, and Proposition \ref{prop: arakelov jensen}. This shows that if $\phi: \Pberk \to \R$ is defined by 
\[\phi(y) := \sum_{\substack{f(\zeta) = \zeta_\Gauss\\ f(\beta) = y}} m_f(\zeta)m_f(\beta) \log \frac{1}{||\zeta, \beta||_v},\] 
then for all $y \in \C_v$
\begin{equation}\label{eq: non arch bounds first step}
    \phi(y) = \psi_y(\zeta_\Gauss) = -\log|a_p^{d-q}b_q^{d-p}\Res(P,Q)|_v + d\log|f|_v + \log^+ |y|_v + d I_v(y).
\end{equation}
Now $f(\zeta) = \zeta_\Gauss$ implies that $\zeta$ is not in $\C_v$. As the spherical kernel is continuous in each variable separately, and since $||\zeta, y||_v \geq ||\zeta, \zeta||_v > 0$ whenever $\zeta$ is not in $\C_v$, then $\phi$ is continuous on $\Pberk$. Suppose that $y_n \in \C_v$ with $y_n \to \zeta_\Gauss$. Applying (\ref{eq: non arch bounds first step}), Proposition \ref{prop: non arch integral limit}, and the fact that $\lim_{y_n \to \zeta_\Gauss} \log^+ |y_n|_v = 0$ yields 
\[\phi(\zeta_\Gauss) = \lim_{y_n \to \zeta_\Gauss} \phi(y_n) = -\log|a_p^{d-q}b_q^{d-p}\Res(P,Q)|_v + 2d \log |f|_v.\]
On the other hand,
\[\phi(\zeta_\Gauss) = \sum_{\substack{ f(\zeta) = \zeta_\Gauss\\ f(\zeta') = \zeta_\Gauss}}  m_f(\zeta)m_f(\zeta')\log\frac{1}{||\zeta, \zeta'||_v} = \|f^* \delta_{\zeta_\Gauss}\|_v^2.\]
\end{proof}

\subsection{Archimedean Considerations}

The Archimedean case is simpler than the non-Archimedean, as we may simply integrate the expression in Lemma \ref{lem: resultant pull back} with respect to $\lambda_{\Ar, v}$. 

\begin{prop}\label{prop: arch pull back of sigma}
At Archimedean places $v$, we have
\[  \|f^* {\lambda_{\Ar, v}}\|_v^2 = -\log|a_p^{d-q}b_q^{d-p}\Res(P,Q)|_v - \frac{d}{2} + 2d \int \log\sqrt{ |P(z)|_v^2 + |Q(z)|_v^2} \; d{\lambda_{\Ar, v}}(z)\]
\end{prop}

\begin{proof}
With notation as in Lemma \ref{lem: resultant pull back},
\begin{align*}
    \|f^* {\lambda_{\Ar, v}}\|^2 &= \iint \sum_{i,j=1}^d \log\frac{1}{||\alpha_i, \beta_j||_v} d{\lambda_{\Ar, v}}(x)d{\lambda_{\Ar, v}}(y)\\
                        &= -\log |a_p^{d-q}b_q^{d-p}\Res(P,Q)|_v - d\iint \log |x-y|_v d{\lambda_{\Ar, v}}(x) d{\lambda_{\Ar, v}}(y) + 2d \int I_v(x)\; d{\lambda_{\Ar, v}}(x).
\end{align*}
(As ${\lambda_{\Ar, v}}$ is log-continuous,  we may safely ignore both the diagonal of $\P^1(\C)\times \P^1(\C)$ and the finite number of points where Lemma \ref{lem: resultant pull back} fails to apply.) The first integral is easily computed:
\[ \iint \log |x-y|_v\; d{\lambda_{\Ar, v}}(x)d{\lambda_{\Ar, v}}(y) = \int \log\sqrt{1 + |y|_v^2} \; d{\lambda_{\Ar, v}}(y) = \frac{1}{2}. \]
For the second:
\begin{align*}
    \int I_v(x) \; d{\lambda_{\Ar, v}}(x) &= \iint \log| P(z) - xQ(z)|_v\; d{\lambda_{\Ar, v}}(z) d{\lambda_{\Ar, v}}(x) \\
                              &= \iint \log |P(z) - xQ(z)|_v\; d{\lambda_{\Ar, v}}(x)\;d{\lambda_{\Ar, v}}(z)\\
                              &= \int \log\sqrt{|P(z)|_v^2 + |Q(z)|_v^2} \; d{\lambda_{\Ar, v}}(z),
\end{align*}  
using Fubini's Theorem in the second line and Proposition \ref{prop: arakelov jensen} in the final line. 
\end{proof}

As one might suspect, it is not difficult to bound the integral $\int \log \sqrt{ |P(z)|_v^2 + |Q(z)|_v^2} \; d{\lambda_{\Ar, v}}(z)$ in terms of the Archimedean contribution to the Arakelov height of $f$. Recall that for Archimedean $v$, we may write $d\lambda_{\Ar, v}(z) = d\ell(z)/(\pi(1 + |z|_v^2)^2$, where $\ell$ is the Lebesgue measure on $\C$.

\begin{prop}\label{prop: bounds on archimedean pull back}
At Archimedean places $v$, we have
\[ \int \sqrt{ |P(z)|_v^2 + |Q(z)|^2_v} \, d\lambda_{\Ar, v}(z) = \log |f|_v + O(1)  \]
where $|f|_v^2 = \sum_{i=0}^d |a_i|_v^2 + |b_i|_v^2$ and where the implied constants depend only on the degree $d$.
\end{prop}

\begin{proof}
Set $\psi(z) = \sqrt{ |P(z)|_v^2 + |Q(z)|_v^2}$. We begin with the upper bound, which is the easier of the two.  Expanding $|P(z)|_v^2 = P(z)\overline{P(z)}$ and $|Q(z)|_v^2= Q(z)\overline{Q(z)}$ shows that
\begin{equation}\label{eq: bounds on arch pull back unit circle} \int_0^{2\pi} \psi^2(re^{it})\; \frac{dt}{2\pi} = \sum_{k=0}^d \left(|a_k|_v^2 + |b_k|_v^2\right)r^{2k}. \end{equation}
Writing $\lambda_{\Ar, v}$ in polar coordinates gives
\begin{align*}    
    \int \log \psi(z) d\lambda_{\Ar, v}(z) &= \int_0^\infty \int_0^{2\pi} \log \psi^2(r e^{it})\frac{dt}{2\pi} \frac{r dr}{(1 + r^2)^2} \\
                                   &\leq \int_0^\infty \log \left( \int_0^{2\pi} \psi^2(re^{it}) \frac{dt}{2\pi} \right) \frac{r dr}{(1 + r^2)^2}  \\
                                   &= \int_0^\infty \log \left( \sum_{k=0}^d \left(|a_k|_v^2 + |b_k|_v^2\right)r^{2k} \right) \frac{r dr}{(1 + r^2)^2}\\
                                   &\leq \int_0^\infty \log \left(|f|_v^2 \sum_{k=0}^d r^2k\right)\frac{r dr}{(1 + r^2)^2}\\
                                   &= \frac{1}{2}\log|f|_v^2 + \int_0^\infty \log \left(\sum_{k=0}^d r^{2k}\right) \; \frac{r\,dr}{(1 + r^2)^2}\\
                                   &=: \log |f|_v + J_d
\end{align*}
where we have used Jensen's inequality in the second line, (\ref{eq: bounds on arch pull back unit circle}) in the third line, and the fact that $\int_0^\infty \frac{r\,dr}{(1+r^2)^2} = \frac{1}{2}$ in the penultimate line. Combining elementary calculus with trivial bounds on $\sum_{k=0}^d r^{2k}$ shows that $J_d \leq \frac{2d + 1}{4}\log(2) + \frac{1}{4}\log(d+1)$. For the lower bound, we assume that $\deg Q = d$ and define a polynomial $H(z)$ by
\[  H(z) = P(z) + z^{d+1}Q(z) = a_0 + a_1z + \cdots + a_d z^d + b_0 z^{d+1} + \cdots + b_d z^{2d+ 1}. \]
Then, by the triangle inequality and the fact that $(A + B)^2 \leq 2A^2 + 2B^2$ for positive numbers $A$ and $B$, we find that $|H(z)|_v^2 \leq 2\max(1, |z|_v^{2d + 2})\psi^2(z)$ for all $z$. Thus

\begin{equation}\label{eq: arch lower bound 1}
\int \log |H(z)|_v\; d\lambda_{\Ar, v} (z) \leq \int_{|z|_v \leq 1} \log \left(\sqrt{2} \psi(z)\right)\; d\lambda_{\Ar, v}(z) + \int_{|z|_v>1} \log\left( \sqrt{2}|z|_v^{d+1}\psi(z)\right) \; d\lambda_{\Ar, v}(z). 
\end{equation}
Now write $H(z) = b_d \prod_{i=1}^{2d+ 1}(z - \gamma_i)$. By Proposition \ref{prop: arakelov jensen}  
\[\int \log |H(z)|_v\; d\lambda_{\Ar, v}(z) = \log |b_d|_v + \sum_{i=1}^{2d+1} \log \sqrt{1 + |\gamma_i|_v^2},\] 
and the roots of $H$ are easy to relate to the coefficients of $H$ via Vieta's formulas: 
\begin{align*}
    |b_d|_v^2\prod_{i=1}^{2d+1} (1 + |\gamma_i|_v^2) &\geq |b_d|_v^2 2^{-(2d+1)}\prod_{i=1}^{2d +1}(1 + |\gamma_i|_v)^2\\
                                                &= 2^{-(2d+1)}\left(  |b_d|_v + |b_d|_v(|\gamma_1|_v + \cdots + |\gamma_{2d+1}|_v) + \cdots + |b_d|_v\prod_{i=1}^{2d+1}|\gamma_i|_v  \right)^2 \\
                                                &\geq 2^{-(2d + 1)}\left( \sum_{i=0}^d |a_i|_v + \sum_{i=0}^d |b_i|_v\right)^2 \\
                                                &\geq 2^{-(2d + 1)}\left(\sum_{i=0}^d |a_i|_v^2 + |b_i|_v^2\right),
\end{align*}
 Consequently,
\begin{equation}\label{eq: arch lower bound 2}  
\int \log |H(z)|_v\; d\lambda_{\Ar, v}(z) = \log |b_d|_v + \sum_{i=1}^{2d+1} \log \sqrt{1 + |\gamma_i|_v^2} \geq -(2d + 1)\log \sqrt{2} + \log |f|_v.  
\end{equation}
From (\ref{eq: arch lower bound 1}), (\ref{eq: arch lower bound 2}), and the fact  $\int_{|z|_v>1} \log |z|_v \; d\lambda_{\Ar, v}(z) = \log \sqrt{2}$, we find  \[\int \log \psi(z) \; d\sigma(z) \geq   - 3(d + 1)\log\sqrt{2} + \log |f|_v.\] 
If instead $\deg Q < d$, then $\deg P = d$ and we set $H(z) = Q(z) + z^{d+1}P(z)$. The above argument then gives the same lower bound on $\int\log \psi\; d\lambda_{\Ar, v}$. 

\end{proof}



Multiplying the estimates in Lemma \ref{lem: non-arch pull back of sigma} and Propositions \ref{prop: arch pull back of sigma} and \ref{prop: bounds on archimedean pull back} by $r_v$, summing over all places and applying the product formula shows that $\frac{1}{2d}\| f^*\lambda_\Ar\|^2 = \frac{d}{2}\|d^{-1}f^* \lambda_\Ar\|^2 = h_{\Ar}(f) + O(1)$, as claimed at the beginning of Section \ref{sec: bounds on pull back}. In fact, the estimates above show that
\begin{equation}\label{eq: explicit arakelov pull back}  - \frac{3(d+1)}{2} \log(2) - \frac{1}{4} + h_\Ar(f) \leq \frac{d}{2}\|d^{-1}f^*\lambda_\Ar\|^2 \leq h_\Ar(f) -\frac{1}{4} + \frac{2d+1}{4} \log(2) + \frac{1}{4}\log(d+1) 
\end{equation}

Combining the estimates in Lemma \ref{lem: pull back bounds} and (\ref{eq: explicit arakelov pull back}) gives explicit bounds on $\|\mu_f\|^2$:
\begin{equation}\label{eq: explicit bounds arakelov height}
    \frac{d}{(\sqrt{d} + 1)^2}h_{\Ar}(f) - c_2  \leq \frac{d}{2}\|\mu_f\|^2 \leq \frac{d}{(\sqrt{d} - 1)^2}h_{\Ar}(f) + c_4
\end{equation}
where $c_2 = \frac{d}{(\sqrt{d} + 1)^2}\cdot(- \frac{3(d+1)}{2} \log(2) - \frac{1}{4} )$ and where $c_4 = \frac{d}{(\sqrt{d}-1)^2} \cdot( -\frac{1}{4} + \frac{2d+1}{4} \log(2) + \frac{1}{4}\log(d+1) )$. This completes the proof of Theorem \ref{thm: arakelov height norm}.





\subsection{The Arakelov--Zhang Pairing is (essentially) a height}
We apply Theorem \ref{thm: arakelov height norm} to prove Theorem \ref{thm: height on polynomials} and Corollary \ref{cor: AZ Pairing is a height} from the introduction.

\begin{cor}\label{cor: ag norm a height}
Let $\lambda = (\lambda_v)_{v \in M_K}$ be an adelic probability measure and let $t > 0$. Then, for any rational map $f$ defined over $K$ of degree $d \geq2$, there are constants $c_i$, $i = 1, 2, 3, 4$ depending on $d$, $\lambda$, and $t$ so that 
\[ c_1 h_{\Ar}(f) - c_2 \leq \frac{d}{2}\|\mu_f\|_{\lambda, t}^2 \leq c_3 h_\Ar (f) + c_4. \]
Moreover $c_1, c_2 \to 1$ as $d \to \infty$ and $c_2, c_4$ grow linearly with $d$. 
\end{cor}

\begin{proof}
    Combine Theorem \ref{thm: arakelov height norm} with the fact that $\|\mu_f\|^2_{\lambda, t} = \|\mu_f\|^2 + O(1)$, where the implied constants depend only on $\lambda$ and $t$ (Proposition \ref{prop: ag inner product essentially standard inner product}).
\end{proof}

\begin{cor}\label{cor: AZ pairing is essentially a height}
    Fix a rational map $g$ with $\deg (g) \geq 2$, and let $f$ be a rational map of degree $d \geq 2$. There are positive constants $c_5, c_6, c_7, c_8$ depending on $g$ and $d$ so that
    \[  c_5 h_{\Ar}(f) - c_6 \leq  d\cdot(f, g)_{\AZ} \leq c_7 h_{\Ar}(f) + c_8.\]
Moreover, $c_5, c_7$ depend only on $d$ and satisfy $c_5, c_7 \to 1$ as $d \to \infty$, and $c_6, c_8$ grow linearly in $d$. If $f$ is a monic polynomial then $d(f, g)_{\AZ} = h_{\Ar}(f) + O(1)$ where the implied constants depend on $g$ and $d$. 
\end{cor}

\begin{proof}
By (\ref{eq: ag norm az pairing}), $(f, g)_{\AZ} = \lim_{t \to 0} \frac{1}{2}\|\mu_f\|_{\mu_g, t}^2$. By Proposition \ref{prop: ag inner product essentially standard inner product} there is a $\kappa$ depending only on $\mu_g$ so that $\big| \|\mu_f\|_{\mu_g, t} - \|\mu_f\|^2 \big|\leq \kappa + t$. Therefore $|2(f, g) - \|\mu_f\|^2| \leq \kappa$, where $\kappa$ depends only on $\mu_g$. Now apply Theorem \ref{thm: arakelov height norm} to $\|\mu_f\|^2$ when $f$ is a rational map, and Theorem \ref{thm: norm of monic poly is a height} when $f$ is a monic polynomial.
\end{proof}

As any Weil height on $\Rat_d(K)$ differs from the Arakelov height by a bounded function, Corollaries \ref{cor: ag norm a height} and \ref{cor: AZ pairing is essentially a height} imply Theorem \ref{thm: height on polynomials} and Corollary \ref{cor: AZ Pairing is a height} from the introduction. 

\subsection{A lower bound for adelic measures with points of small height}
In this section, we prove the statements in Remark \ref{rem: small points intro} from the introduction. What follows summarizes some of the discussion in \cite{favre2006equidistribution} and \cite{fili2017metric}. Given an element $\alpha \in \P^1(\bar{K})$ and a place $v\in M_K$,  we denote by $[\alpha]_v$ the probability measure on $\Pberk$ supported equally on the Galois conjugates of $\alpha$ over $K$; thus 

\[ [\alpha]_v = \frac{1}{|\Gal(\bar{K}/K)\cdot \alpha| }\sum_{x \in \Gal(\bar{K}/K)\cdot \alpha} \delta_x.\]

Each $x \in \Gal(\bar{K}/K)\cdot \alpha$ is viewed as an element of $\Pberk$ via an embedding $\bar{K} \hookrightarrow \C_v$; averaging over the Galois conjugates of $\alpha$ ensures that the choice of embedding is irrelevant. We write $[\alpha] = ([\alpha]_v)_{v \in M_K}$. Recall that the mutual energy pairing of two adelic measures $\mu = (\mu_v)_{v \in M_K}, \nu = (\nu_v)_{v \in M_K}$ is
\[ (\!( \mu, \nu) = \sum_{v \in M_K} r_v (\!( \mu_v, \nu_v)\!)_v, \; \text{where}\; (\!(\mu_v, \nu_v)\!)_v = \iint_{\Aberk\times \Aberk\setminus\Diag_v} \log \delta_v(x,y)^{-1} \;d\mu_v(x)\; d\nu_v(y). \]

Recall that $\mu$ \textit{has points of small height} if there is a sequence $\{\alpha_n\}_{n=1}^\infty$ of distinct points in $\P^1(\bar{K})$ so that $h_{\mu}(\alpha_n) \to 0 \;\;\text{as}\;\; n \to \infty$, where $h_{\mu}(\alpha) := \frac{1}{2} (\!( \mu - [\alpha], \mu - [\alpha])\!)$ is the Favre--Rivera-Letelier height associated to $\mu$. The \textit{Arakelov height} of a point $\alpha  \in \P^1(\bar{K})\setminus\{\infty\}$ is defined to be 
\[ h_{\Ar}(\alpha) = \frac{1}{n}\sum_{i=1}^n \sum_{v \in M_K} r_v \log ||\alpha_i, \infty||_v^{-1} \]
where $\alpha_1, \alpha_2, ..., \alpha_n$ are the $\Gal(\bar{K}/K)$-conjugates of $\alpha$. For $\alpha = \infty$, we set $h_\Ar(\infty) =0$. The Arakelov height and the Arakelov measure are related by 
\begin{equation}\label{eq: height of arakelov measure and arakelov height}
    h_{\lambda_\Ar}(\alpha) = h_{\Ar}(\alpha) - \frac{1}{4}
\end{equation}
for all $\alpha \in \P^1(\bar{K})$. This is because
    \[ (\!( \lambda_\Ar, \lambda_\Ar)\!) = \frac{-1}{2}, \;\; -(\!( \alpha, \lambda_\Ar)\!) = h_{\Ar}(\alpha), \;\; \text{and}\;\; (\!([\alpha], [\alpha])\!) = 0.\] 
The first equality is a consequence of Proposition \ref{prop: arakelov jensen}, the second follows from Proposition \ref{prop: arakelov jensen} and the definition of the Arakelov height, and the third is by the product formula. The \textit{standard adelic measure}, $\lambda_\std = (\lambda_{\std, v})_{v \in M_K}$ is defined to be the Haar measure on the unit circle when $v$ is Archimedean, and $\lambda_{\std, v} = \delta_{\zeta_\Gauss}$ when $v$ is non-Archimedean. A direct calculation shows that $\lambda_\std$ is invariant under the pull-back by a power map. Therefore $\lambda_\std$ is the canonical adelic measure of $z \mapsto z^d$ for $|d| \geq 2$. 

\begin{prop}\label{prop: small points}
Let $\mu$ be an adelic probability measure and let $\{\alpha_n\}_{n=1}^\infty$ a sequence in $\P^1(\bar{K})$ with $h_\mu(\alpha_n) \to 0$. Then $h_\Ar(\alpha_n) \to \frac{1}{2}\|\mu\|^2$ as $n \to \infty$. Moreover $\|\mu\|^2 \geq \log(2)$ with equality if and only if $\mu = \lambda_\std$. 
\end{prop}

\begin{proof}
    Suppose that $\{\alpha_n\}_{n=1}^\infty$ is a sequence of distinct points in $\P^1(\bar{K})$ so that $h_{\mu}(\alpha_n) \to 0$. Then Theorem 9 of \cite{fili2017metric} implies $h_{\lambda_\Ar}(\alpha_n) \to \frac{1}{2}\|\mu- \lambda_\Ar\|^2 = \frac{1}{2}\|\mu\|^2 - \frac{1}{4}$. So (\ref{eq: height of arakelov measure and arakelov height}) implies $h_{\Ar}(\alpha_n) \to \frac{1}{2}\|\mu\|^2$ as $n \to \infty$. \\
    Now $h_\Ar(\alpha) \geq \frac{1}{2}\log(2)$ for all $\alpha \in \P^1(\bar{K}) \setminus\{0, \infty\}$; see \cite{fili2017energy}, \cite{abarcelona2005minimums}. Therefore $\frac{1}{2}\|\mu\|^2 = \lim_{n \to \infty} h_{\Ar}(\alpha_n) \geq \frac{1}{2}\log(2)$. A straightforward computation shows that $\|\lambda_\std\|^2= \log(2)$; see Proposition \ref{prop: norm of conjugation of squaring map} for a more general statement. Lastly, if $\|\mu\|^2 = \log(2)$ then $h_{\Ar}(\alpha_n) \to \frac{1}{2}\log(2)$. It is known that this implies $[\alpha_n]_v \to \lambda_{\std, v}$ weakly at each place of $K$; see Example 6.5 of \cite{gil2019distribution}. On the other hand, Théorème 2 of \cite{favre2006equidistribution} implies $[\alpha_n]_v$ converges weakly to $\mu_v$ at each place. So $\mu = \lambda_\std$. 
\end{proof}

\section{Examples}\label{sec: examples}

\subsection{Chebyshev Polynomials}

Let $T_n$ be the $n^{th}$ Chebyshev polynomial. We adopt the convention that $T_2(z) = z^2 - 2$ so that $T_2$ (and consequently $T_n$) has good reduction at all finite places. In general, the (Archimedean) canonical measure of a polynomial is the equilbrium measure of the Julia set, a result due to Brolin \cite{brolin1965invariant}. It is well known that the complex Julia set of $T_n$ is the interval $[-2, 2]$. The map $\phi(z) = z+ z^{-1}$ is a conformal mapping from $\P^1(\C) \setminus \{ |z| \leq 1\}$ onto $\P^1(\C) \setminus [-2, 2]$. By potential theory, this implies that $\mu_{T_n, \infty}$ is the push forward by $\phi$ of the Haar measure on the unit circle; see Theorem 4.4.4 of \cite{ransford1995potential} for an equivalent statement using Green's functions. With this explicit description for the canonical measure of $T_n$, it is possible to compute the norm of $\mu_{T_n}$ exactly. 
\begin{prop}
Let $\mu_{T_n}$ be the canonical adelic measure of $T_n$. Then  $\|\mu_{T_n}\|^2 = \log \frac{3 + \sqrt{5}}{2}$.
\end{prop}

\begin{proof}
Proposition \ref{prop: potential at infinity} implies that $\|\mu_{T_n}\|^2 = 2\sum_{v \in M_\Q}  U_{T_n, v}(\infty)$. Since $T_n$ has good reduction at all finite places, then $\mu_{T_n, v} = \delta_{\zeta_\Gauss}$ for $v \nmid \infty$; consequently $U_{T_n,v}(\infty) =0$ when $v \nmid \infty$. So
\[ \|\mu_{T_n}\|^2 = 2U_{T_n, \infty}(\infty)  =  \int \log(1 + |z|^2) \; d\mu_{T_n, \infty}(z) = \int_0^{2\pi} \log \left(1 + |e^{it} + e^{-it}|^2 \right) \; \dfrac{dt}{2\pi}.  \]
Now $1 + |e^{it} + e^{-it}|^2 = |re^{it} + r^{-1}e^{-it}|^2$, where $r^2 = \frac{3 + \sqrt{5}}{2}$. Therefore
\[ \|\mu_{T_n}\|^2 = \int_0^{2\pi} \log |re^{-it}|^2 | e^{2it}+  r^{-2}|^2\;\dfrac{dt}{2\pi} = \log r^2 + \int_0^{2\pi} \log |e^{2it} + r^{-2}|^2\; \dfrac{dt}{2\pi}.  \]
By Jensen's formula, this final integral vanishes. So $\|\mu_{T_n}\|^2 = \log r^2$. 

\end{proof}

\subsection{The Squaring Map After an Affine Change of Variables} 
In this section we compute the norm of the canonical measure associated to the map $P_d^\phi = \phi^{-1} \circ P_d \circ \phi$, where $P_d: z \mapsto z^d$ is the $d^{th}$ power map ($|d|\ \geq 2$) and where $\phi : z \mapsto az + b$ is an affine change of variables defined over $K$.  In general, the canonical measure of a rational map $f$ of degree $d$ after conjugation by $\phi\in \PGL_2(K)$ can be found via the formula

\begin{equation}\label{eq: conjugation formula} 
\phi^* \mu_{f,v} = \mu_{f^\phi,v}. 
\end{equation}

To see this, note that $(f^\phi)^* = \phi^* f^* (\phi^{-1})^*$, and also that $(\phi^{-1})^* = \phi_*$, as $\phi$ is invertible. So 
\[ (f^\phi)^* (\phi^* \mu_{f,v}) = \phi^* f^* (\phi_*\phi^* \mu_{f, v}) = \phi^* f^*\mu_f = d\phi^*\mu_{f,v}. \]
This and the fact that the canonical measure is characterized by the identity $(f^\phi)^*\mu_{f^\phi, v} = d \mu_{f^\phi}$ shows (\ref{eq: conjugation formula}). As remarked above, the canonical adelic measure of a power map is $\lambda_\std$; this adelic measure is simple enough to admit a direct computation of $\|\mu_{P_d^\phi}\|^2 = \|\phi^*\lambda_\std\|^2$.

\begin{prop}\label{prop: norm of conjugation of squaring map}
Let $\phi(z) = az+ b$ with $a \in K\setminus\{0\}$ and $b \in K$ and let $P_d(z) = z^d$ $\;(|d| \geq 2)$. Then 

\[ \|\mu_{P_d^\phi}\|^2 = -\log 2 +  \sum_{v \in M_K^0} r_v \log\max(1, |a|^2_v, |b|^2_v) + \sum_{v \in M_K^\infty} r_v\log \eta_v (a, b),\]

where $\eta_v (a, b) = 1 + |a|_v^2 + |b|_v^2 + \sqrt{( |a|_v^2 + (1 - |b|_v)^2)( |a|_v^2+ (1 + |b|_v)^2)}$.
\end{prop}

\begin{proof}
From Proposition \ref{prop: potential at infinity}, $2\sum_{v \in M_K} r_v U_{\mu, v}(\infty) = \|\mu_{P_d^\phi}\|^2$. So we compute $2U_{P_d^\phi, v}(\infty)$ at each place. For non-Archimedean $v$, (\ref{eq: conjugation formula}) implies
\[ U_{{P_d^\phi}, v}(\infty) = \int \log^+|z|_v \; d\phi^*\delta_{\zeta_\Gauss} = \int \log^+ | \phi^{-1}(z)|_v \; d\delta_{\zeta_\Gauss}(z),  \]
where, as in Section \ref{sec: comparison w ar ht rat}, $|\cdot|_v$ is the continuous extension of the absolute value from $\C_v$ to $\Aberk$. Thus $|\phi^{-1}(\zeta_\Gauss)|_v = \sup_{\substack{z \in \C_v \\|z|_v \leq 1}} |(z - b)/a|_v = \max( |a|_v^{-1}, |b/a|_v)$ and so 

\begin{equation}\label{eq: conjugation square non arch}
    2U_{{P_d^\phi}, v}(\infty) = \log \max( 1, |a|_v^{-2}, |b/a|^2_v).
\end{equation}

For Archimedean $v$, we are led to consider integrals of the form $\int_0^{2\pi} \log(1 + |\alpha e^{i\theta} + \beta|^2)  \frac{d\theta}{2\pi}$, 
where $\alpha$ and $\beta$ are complex numbers. We claim that

\begin{equation}\label{eq: conjugation square arch}
\int_0^{2\pi} \log \left( 1 + |\alpha e^{i\theta} + \beta|^2\right)\; \frac{d\theta}{2\pi} = \log \left(1 + \alpha^2 + \beta^2 + \sqrt{(1 + (\alpha - \beta)^2)(1 + (\alpha + \beta)^2)}\right) - \log 2.
\end{equation}
Assuming (\ref{eq: conjugation square arch}) is true, then  (\ref{eq: conjugation formula}) implies

\begin{align*}
    2U_{{P_d^\phi}, v}(\infty) &= \int_{0}^{2\pi} \log(1  + |\phi^{-1}(e^{i\theta})|_v^2) \; \frac{d\theta}{2\pi} \\
                               &= \int_0^{2\pi} \log( 1 + |(e^{i\theta} - b)/a|_v^2 ) \frac{d\theta}{2\pi}\\
                                 &= \log \left(1 + |a|_v^{-2} + |b/a|_v^2 + \sqrt{(1 + (|a|_v^{-1} - |b/a|_v )^2)(1 + (|a|_v^{-1} + |b/a|_v)^2)}\right) - \log 2.\\
\end{align*} 
This, (\ref{eq: conjugation square non arch}), and the product formula imply
\begin{align*}
    \| \mu_{P_d^\phi}\|^2 &= 2\sum_{v\in M_K} r_v U_{{P_d^\phi}, v}(\infty)\\
                          &= \sum_{v\in M_K}r_v \log |a|_v^2 + 2\sum_{v\in M_K} r_v U_{{P_d^\phi}, v}(\infty)\\
                          &= \sum_{v \in M_K^0} r_v \log\max(1, |a|^2_v, |b|^2_v) + \sum_{v \in M_K^\infty} r_v(\log \eta_v (a, b) - \log 2),
\end{align*}
where $\eta_v(a,b) = 1 + |a|_v^2 + |b|_v^2 + \sqrt{(|a|_v^2 + (1 - |b|_v)^2)(|a|_v^2+ (1 + |b|_v)^2)}$. Applying the fact that $\sum_{v \in M_K^\infty}r_v (-\log 2) = -\log 2$ yields the proposition. We therefore turn to proving (\ref{eq: conjugation square arch}). Note  that the integral in (\ref{eq: conjugation square arch}) is unchanged if $\alpha$ and $\beta$ are replaced by $\alpha e^{i\theta'}$ and $\beta e^{i\theta''}$ with $\theta', \theta''\in [0, 2\pi]$. So assume without loss of generality that $\alpha = |\alpha|$ and $\beta = |\beta|$. It is convenient to define the quantities
\[ A = \sqrt{ \frac{1 + (\alpha - \beta)^2}{2}} \qand B = \sqrt{ \frac{1 + (\alpha + \beta)^2}{2}}, \]
and 
\[ r = \frac{1}{2}\left(A + B\right)^2 \qand s = \frac{B^2 - A^2}{(A + B)^2}.  \]
We have the following identities: 
\[ 1 + \alpha^2 + \beta^2 = A^2 + B^2, \;\; 2\alpha\beta = B^2 - A^2, \qand 1 + s^2  = 2\frac{A^2 + B^2}{(A + B)^2} . \]
It follows from these equations that $1+ |\alpha e^{i\theta} + b|^2  = r|1 + se^{i\theta}|^2$ for all $\theta \in [0, 2\pi]$. Consequently

\[ 
\int_0^{2\pi}\log(1 + |\alpha e^{i\theta} + \beta|^2)\;\frac{d\theta}{2\pi} = \log r + \int_0^{2\pi} \log|1 + se^{i\theta}|^2\; \frac{d\theta}{2\pi}. \]
This final integral is zero by Jensen's formula: 
\[ \int_0^{2\pi} \log|1 + se^{i\theta}|^2\; \frac{d\theta}{2\pi} = \log|s|^2 + \log^+|1/s|^2 = \log^+ |s|^2= 0.  \]
So the integral in (\ref{eq: conjugation square arch}) is equal to $\log r$. Expressing $r$ in terms of $\alpha$ and $\beta$ gives the claim.
\end{proof}

\subsection{A Family of Latt\`es Maps} 
Following \cite{petsche2012dynamical}, we consider the elliptic curves
\[ E_{a, b}: \;  y^2 = P(x) =  x(x- a)(x+b), \]
where $a$ and $b$ are positive integers. The Latt\`es map associated to $E_{a, b}$ is then 
\begin{equation}\label{eq: lattes map def} 
f_{a, b}(x) = \frac{ (x^2 + ab)^2  }{4x(x-a)(x+b)}. 
\end{equation}
This is the  rational map $f_{a,b} : \P^1 \to \P^1$ which is semi-conjugate via the $x$-coordinate function to the doubling map on $E_{a, b}$. That is,  if $P$ is a point on $E_{a, b}$ and $[2]: E_{a, b} \to E_{a, b}$ is the duplication map, then $f_{a, b}(x(P)) = x( [2]P)$.  Proposition 20 of \cite{petsche2012dynamical} states that

\[ \frac{1}{2}\log(ab) \leq (f_{a, b}, z^2 )_{\AZ} \leq  \frac{1}{2}\log(ab) + C \]
where $C$ is some absolute positive constant. Proposition \ref{prop: ag inner product essentially standard inner product} and (\ref{eq: ag norm az pairing}) imply that $(f_{a,b}, z^2)_{\AZ} = \frac{1}{2}\|\mu_{f_{a,b}}\|^2 + O(1)$. Therefore
\begin{equation}\label{eq: lattes o one}
\frac{1}{2} \|\mu_{f_{a,b}}\|^2 = \frac{1}{2}\log(ab) + O(1)
\end{equation}
An inspection of (\ref{eq: lattes map def}) reveals that the naive Weil height of the coefficients of $f_{a, b}$ is $\log(a^2 b^2)$ for $a$ and $b$ sufficiently large (for instance, when $ab \geq 4$). So (\ref{eq: lattes o one}) can be rephrased as $2\|\mu_{f_{a, b}}\|^2 = h(f_{a, b}) + O(1)$. This provides an example of an infinite family of rational maps for which $\frac{\deg(f)}{2}\|\mu_f\|^2 = h(f) + O(1)$, which is a stronger statement than that of Theorem \ref{thm: height on polynomials}. This suggests that it is possible that Theorem \ref{thm: height on polynomials} could be strengthened to show the estimate $\frac{\deg(f)}{2}\|\mu_f\|^2 = h(f) + O(1)$ holds for rational maps, and not just monic polynomials.

\printbibliography

\end{document}